\documentclass[12pt,reqno]{article}
\usepackage{amssymb, amsthm, amsmath, amsfonts, stmaryrd, latexsym, txfonts, pxfonts, wasysym}
\usepackage[v2]{xy}
\xyoption{curve}

\usepackage[usenames]{color}
\usepackage{amssymb}
\usepackage{graphicx}
\usepackage{amscd}

\usepackage{fullpage}
\usepackage{setspace}
\usepackage[colorlinks=true,
linkcolor=webgreen,
filecolor=webbrown,
citecolor=webgreen]{hyperref}

\definecolor{webgreen}{rgb}{0,.5,0}
\definecolor{webbrown}{rgb}{.6,0,0}

\usepackage{color}
\usepackage{fullpage}
\usepackage{float}

\usepackage{graphicx,amsmath,amssymb}
\usepackage{amsthm}
\usepackage{amsfonts}
\usepackage{latexsym}
\usepackage{epsf}

\newcommand{\al}{\alpha}

\newcommand{\ga}{\gamma}
\newcommand{\de}{\delta}

\newcommand{\Pp}{{\mathcal P}}

\newcommand{\iso}{\cong}   		% preferred isomorphism symbol
\newcommand{\id}{\text{id}}

\newcommand{\Cc}{{\mathcal C}}
\newcommand{\HC}{\text{Hom}_{\mathcal C}}
\newcommand{\HCT}{\text{Hom}_{{\mathcal C}(\tau)}}
\newcommand{\dm}{\sslash}

\begin{document}

\begin{center}
\vskip 1cm{\LARGE\bf
A new semidirect product for association schemes}

\vskip 1cm
\large
Christopher French\\
Department of Mathematics and Statistics\\
Grinnell College\\
Grinnell, IA 50112\\
USA\\
\href{mailto:frenchc@grinnell.edu}{\tt frenchc@grinnell.edu}\\
\end{center}

\vskip .2in

\begin{abstract}
We propose a new definition for the semidirect product
of association schemes, generalizing work of Bang, Hirasaka, and Song.
We define an action of one scheme on another, and show how one can use such an action
to construct all semidirect products satisfying a certain technical condition.
\end{abstract}

\newtheorem{thm}{Theorem}[section]
\newtheorem{cor}[thm]{Corollary}
\newtheorem{lem}[thm]{Lemma}
\newtheorem{question}[thm]{Question}
\newtheorem{con}[thm]{Conjecture}
\newtheorem{prop}[thm]{Proposition}
\theoremstyle{definition}
\newtheorem{defn}[thm]{Definition}
\newtheorem{notn}[thm]{Notation}
\newtheorem{rem}[thm]{Remark}
\newtheorem{example}[thm]{Example}
\newtheorem{cond}[thm]{Condition}

\section{Introduction}

We begin with an overview of the theory of association schemes, and the connection of that theory
with groups.  Details for what follows can be found in Zieschang's Theory of Association Schemes, \cite{PHZ}.

Given a set $X$, an association scheme on $X$, or a scheme on $X$ for short, is a set $S$ consisting
of nonempty subsets of $X\times X$, satisfying the following axioms:
\begin{enumerate}

	\item $S$ is a partition of $X\times X$
	
	\item The diagonal subset $1_X\subseteq X\times X$, defined as $\{(x,x):x\in X\}$, is an element
	in $S$.
	
	\item For each element $s\in S$, the set $s^*=\{(y,x)\in X\times X:(x,y)\in s\}$ is an element in $S$.
	
	\item For each triple $p, q, r\in S$, and each element $(x,y)\in r$, the cardinality of the set of elements $z\in X$
	such that $(x,z)\in p$ and $(z,y)\in q$ depends only on $p$, $q$, and $r$ (and not on $x$ and $y$).
	The cardinal number for this set is denoted $a_{pqr}$.
	
\end{enumerate}
We will assume throughout that $X$ is a finite set, so that each $a_{pqr}$ is a non-negative integer.

For example, suppose given a finite group $G$.  For each $g\in G$, let $s_g$ be the set of all pairs
$(x,y)\in G\times G$ such that $y=xg$.  Let $S_G=\{s_g:g\in X\}$.  Then

	\begin{enumerate}
	
	\item $S_G$ is a partition of $G\times G$, since for each pair $(x,y)\in G$, we have $(x,y)\in s_g$ if and only
	if $g=x^{-1}y$.
	
	\item $1_G$ is equal to $s_1$, where $1\in G$ is the identity element.
	
	\item For each $g\in G$, $s_g^*=s_{g^{-1}}$, since $y=xg$ if and only if $x=yg^{-1}$.
	
	\item For each triple $g_1, g_2, g_3\in G$, and each pair $(x,y)\in s_{g_3}$, the number of elements
	$z\in G$ such that $(x,z)\in s_{g_1}$ and $(z,y)\in s_{g_2}$ is $1$ if $g_1g_2=g_3$ and $0$
	otherwise.
	
	\end{enumerate}
	
The scheme $S_G$ on $G$ has the property that for
each $p\in S$ and $x\in G$, there is exactly one $y\in G$ such that $(x,y)\in p$.  Schemes with this
property are called thin schemes, and any thin scheme is isomorphic to a scheme
obtained from a group as described above.  Thus, schemes may be
viewed as a generalization of groups.

Many concepts in group theory, like products, subgroups, and quotient groups,
have corresponding generalizations to scheme theory.
One can define the product $pq$ of two elements $p$ and $q$
in a scheme $S$, but in this case, $pq$ is not an element of $S$, but a subset of $S$;
that is $pq$ is the set of all $r$ such that $a_{pqr}>0$.  Similarly,
given subsets $P$ and $Q$ in $S$, we can define $PQ$ to be the union
over all pairs $p\in P, q\in Q$ of the sets $pq$.  We also write $pQ$ for $\cup_{q\in Q} pq$
and $Pq$ for $\cup_{p\in P} pq$.

Corresponding to subgroups of a group are the closed subsets of a scheme.
A subset $T\subseteq S$ is said to be a closed
subset if $T^*T\subseteq T$, where $T^*=\{t^*\in S:t\in T\}$.  If $S_G$ is the thin scheme arising from a group
$G$, then the closed subsets of $S_G$ are precisely the sets of the form $S_H=\{s_g\in S_G:g\in H\}$,
where $H$ is a subgroup of $G$.  As with groups, we say that a closed subset $T$ is normal if $pT=Tp$
for all $p\in S$.

Just as subgroups partition groups into cosets, a closed subset $T$ of a scheme $S$ on a
set $X$ partitions $X$ into a set of cosets.  For $x\in X$ and $t\in T$, we write $xt$
for the set of all $y\in X$ such that $(x,y)\in t$.  We then write $xT$ for $\cup_{t\in T}xt$.
Thus, $x_1$ and $x_2$ are in the same coset if the pair $(x_1,x_2)$ belongs to an element in $T$.
We denote the set of all cosets of $T$ in $X$ as $X/T$.
One obtains a scheme $T_{xT}$ on $xT$, called {\it the subscheme of $S$ defined by $xT$}, whose elements are
the sets $t\cap (xT\times xT)$, where $t\in T$.

For each $s\in S$, we define $s^T$ to be the set of all pairs $(x_1T,x_2T)\in X/T\times X/T$
such that $(x_1',x_2')\in s$ for some $x_1'\in x_1T$ and $x_2'\in x_2T$.  The set of all such
sets $s^T$ then forms a scheme $S\dm  T$ on $X/T$.  In contrast to group theory, a scheme
admits a quotient for any closed subset $T$, not just for normal closed subsets.

Finally, we can define a morphism between two schemes.  Given schemes $S$ on $X$ and $S'$ on $X'$,
a morphism $\phi:S\to S'$ is a function $\phi_X:X\to X'$ such that if $(x_1, x_2)$ and $(x_1', x_2')$ belong
to the same scheme element in $S$, then $(x_1\phi_X,x_2\phi_X)$ and $(x_1'\phi_X,x_2'\phi_X)$ belong to the
same scheme element in $S'$.  Note that $\phi_X$ gives rise to a function $\phi_S$ from the set of scheme elements
$S$ to the set of scheme elements $S'$:  given $s\in S$, choose $(x_1, x_2)\in s$, and
define $s\phi_S$ to be the scheme element in $S'$ containing $(x_1\phi_X, x_2\phi_X)$.  Our definition
ensures that $\phi_S$ is well-defined.  
As an example, if $T$ is a closed subset of a scheme $S$ defined on $X$, then there is a canonical
morphism $\pi:S\to S\dm  T$.  The function $\pi_X:X\to X\dm  T$ takes a given $x\in X$ to its coset $xT$.  The function
$\pi_S$ takes $s$ to $s^T$.  Clearly, if $(x,y)\in s$, then $(xT,yT)\in s^T$.
If $\phi:S\to S'$ is a morphism of schemes, we define the kernel of $\phi$ to be the set of $s\in S$
such that $s\phi=1_{X'}$.  As with groups, there is
an induced morphism $\bar\phi:S\dm T\to S'$, given by $(xT)\bar\phi_{X/T}=x\phi$
and $(s^T)\bar\phi_{S\dm T}=s\phi_S$.  Thus, $\pi \bar\phi=\phi$.
Given a morphism of schemes $\phi:S\to S'$, we will often use $\phi$ to denote
$\phi_X$ or $\phi_S$, relying on context to make the meaning clear.

Since schemes generalize groups, it is natural to look for appropriate generalizations of various
concepts from group theory to scheme theory.  One classical problem in group theory is the extension problem.
Given a pair of groups $N$ and $H$, one seeks to classify those groups $G$
with a given normal subgroup $\tilde N$ isomorphic to $N$, and quotient $G/\tilde N$ isomorphic to $H$.
Such a group is called an extension of $N$ by $H$.  In the same way, one could pose the following problem:
given two schemes $T$ and $U$, classify the schemes $S$ which are equipped with a closed subset $\tilde T$,
such that the quotient scheme $S\dm  \tilde T$ is isomorphic to $U$, and such that the subscheme of
$S$ defined by one of the cosets of $\tilde T$ is isomorphic to the scheme $T$.

The simplest class of extensions of a group $N$ by another group $H$ is the class of
semidirect products of $N$ and $H$.  In considering the extension problem for schemes, it is
natural to seek a generalization of the semidirect product, which should satisfy a certain
collection of properties similar to those of a semidirect product of groups.  Bang, Hirasaka, and
Song \cite{BHS} have already proposed one construction for semidirect products.
Our goal in this paper is to propose a more general construction.

To understand our motivation, we must examine semidirect products of groups more closely.
In group theory, given a pair of groups $H$ and $N$ together with a homomorphism $\zeta:H\to \text{Aut}(N)$,
one defines the semidirect product $G=N\rtimes_\zeta H$ as the group whose underlying set of elements is
the Cartesian product $N\times H$, with product defined by \[(n_1,h_1)(n_2,h_2)=(n_1 (n_2(h_1\zeta)),h_1h_2).\]
Then $G$ contains a normal subgroup $\tilde N=\{(n,1_H):n\in N\}$,
which is isomorphic to $N$, and $G/\tilde N$ is isomorphic to $H$.
We have the following properties:
\begin{enumerate}
	\item $G$ contains a subgroup $\tilde H$ isomorphic to $H$,
	$G=\tilde N\tilde H$ and $\tilde H\cap \tilde N=\{1\}$,
	\item The homomorphism $\pi:G\to G/\tilde N$ splits, in the sense that there is a homomorphism
	$i:H\to G$ such that $i\pi:H\to G/\tilde N$ is an isomorphism.
\end{enumerate}
In fact, given a group $G$ containing a normal subgroup $\tilde N$ isomorphic to $N$,
such that $G/\tilde N$ is isomorphic to $H$, then properties (1) and (2) are equivalent,
and each implies that $G$ is isomorphic to $N\rtimes_\zeta H$ for some $\zeta:H\to \text{Aut}(N)$.

Returning to schemes, there are then two approaches one might take to generalize semidirect products.
On the one hand, one could say that a scheme $S$ on a set $X$ is a semidirect product of $T$ by $U$ if
$S$ contains a normal closed subset $\tilde T$ with a subscheme isomorphic to $T$,
as well as a closed subset $\tilde U$  isomorphic to $U$, such that $\tilde U\cap \tilde T=\{1_X\}$
and $\tilde T\tilde U=S$.
Bang, Hirasaka, and Song \cite{BHS} describe a construction that takes two schemes $T$ and $U$,
together with a kind of twisting map $\zeta$, and outputs a scheme $S=T\rtimes_\zeta U$
with the above properties.  Zieschang \cite[7.4]{PHZ} proves that any scheme satisfying
the conditions above, together with one extra condition (which holds automatically if $T$ is thin),
is isomorphic to a scheme arising from their construction.  The extra required
condition is \begin{cond}\label{con:hypb} For any $u\in \tilde U$ and $t\in \tilde T$, we have $|tu|=1$.
\end{cond}

We now present a second approach to defining semidirect products, which we will adopt in this paper.

\begin{defn}\label{def:semi}
A scheme $S$ on a set $X$ is a {\it semidirect product of $T$ by $U$} if
$S$ contains a closed subset $\tilde T$,
there is a morphism $i:U\to S$ such that $i\pi:U\to S\dm \tilde T$ is an isomorphism,
and the subscheme of $S$ defined by one of the cosets of $\tilde T$
is isomorphic to $T$.
\end{defn}

Note that the image of a scheme morphism $i:U\to S$ need not determine a closed subset of $S$;
that is, the image of $i_{U}:U\to S$ need not be closed.  On the other
hand, if $S$ contains a normal closed subset $\tilde T$, as well as a closed subset $\tilde U$
isomorphic to $U$,
such that $\tilde U\cap \tilde T=\{1_X\}$, and $\tilde T\tilde U=S$, then
the inclusion of the subscheme of $S$ defined by any of the cosets of $\tilde U$, composed
with the projection $\pi:S\to S\dm \tilde T$, is an isomorphism of schemes, and thus determines a splitting
of $\pi$.  Therefore, Definition \ref{def:semi} is more general than
the semidirect product defined by Bang, Hirasaka, and Song.

Our goal in this paper is to show how one can take two schemes $T$ and $U$,
together with a kind of twisting map $\zeta$, and construct an extension $S$ of $T$
by $U$ which is a semidirect product in the sense of Definition \ref{def:semi}.  We will show
that every semidirect product satisfying one extra condition (analogous to 
Condition \ref{con:hypb} above) is isomorphic
to a scheme arising from our construction.  Our extra condition is
\begin{cond}\label{con:hyp} for any $u\in U$ and $t\in \tilde T$, we have $|t(ui)|=1$.
\end{cond}

Bang and Song \cite{BS} also give a generalization of their own semidirect product.  
As we show in Section \ref{sec:example}, it is possible to describe a semidirect product using
our approach that cannot be obtained as a generalized semidirect product in their sense.
On the other hand, any generalized semidirect product in their sense would
arise from our construction.

In Sections \ref{sec:cat} and \ref{sec:label}, we work through some preliminary definitions which
we will need in order to construct our semidirect product.  In
Section \ref{sec:act}, we define an action of one scheme on another, and show how this
can be used to construct a semidirect product, much the way one can use a homomorphism
from a group $G$ to the automorphisms of a group $N$ to construct a semidirect product in
group theory.  We discuss properties of our semidirect product in Section \ref{sec:prop}
and provide an example in Section \ref{sec:example}.  In the final section, Section \ref{sec:action},
we show how, given Condition \ref{con:hyp}, one can obtain an action of one scheme on another
from a given semidirect product, and then reconstruct this given semidirect product
up to isomorphism using our construction.

\section{A category of schemes}\label{sec:cat}

A based set $X$ is a set together with a distinguished element $x_*\in X$, 
called the basepoint.  A based scheme is a scheme $T$ on a based set $X$.
If $T$ and $U$ are based schemes (on based sets $X$ and $Y$), then
a morphism of schemes $\phi:T\to U$ is a based morphism if it takes
the basepoint of $X$ to the basepoint of $Y$.  Any quotient of a based
scheme is naturally a based scheme, where we take the basepoint of the quotient
to be the coset containing the basepoint of the original scheme.

Of course, one could define a category whose objects were based schemes and
whose morphisms were the above morphisms of based schemes.
We now define a different category $\Cc$ of based schemes.  An object of $\Cc$ is a
based scheme $T$ on $X$.
Given two based schemes $T$ and $U$, a morphism $\phi\in \HC (T,U)$ consists of
\begin{enumerate}
	\item a normal closed subset $T_\phi$ in $T$,
	\item a normal closed subset $U_\phi$ in $U$,
	\item a based isomorphism of schemes $\tilde\phi:T\dm  T_\phi\to U\dm U_\phi.$
\end{enumerate}

To define composition of morphisms in $\Cc$, suppose given $\phi\in \HC (T,U)$
and $\psi\in \HC (U,V)$, where $T$, $U$, and $V$ are schemes on $X$, $Y$, and
$Z$ respectively.  We define
\[T_{\phi\psi}=\{t\in T:\text{ for some } u\in U\text{ we have }
t^{T_\phi}\tilde\phi=u^{U_\phi}\text{ and }
u^{U_\psi}\tilde\psi=1_Z^{V_\psi}\},\]
\[V_{\phi\psi}=\{v\in V:\text{ for some } u\in U\text{ we have }
1_X^{T_\phi}\tilde\phi=u^{U_\phi}\text{ and }
u^{u_\psi}\tilde\psi=v^{V_\psi}\}.\]

In Lemma \ref{lem:comp} below, we will show that $T_{\phi\psi}$ and $V_{\phi\psi}$ are closed subsets.
Now, we define a function $\widetilde{\phi\psi}_{X/T_{\phi\psi}}:X/T_{\phi\psi}\to Z/V_{\phi\psi}$
as follows.  Given $xT_{\phi\psi}\in X/T_{\phi\psi}$, choose a coset representative $y$
so that $(xT_{\phi})\tilde\phi=yU_{\phi}$.  Then, choose a coset representative
$z$ so that $(yU_{\psi})\tilde\psi=zV_{\psi}$.  Let
$(xT_{\phi\psi})(\widetilde{\phi\psi}_{X/T_{\phi\psi}})=zV_{\phi\psi}.$
We will let $\phi\psi\in \HC (T,V)$ be the morphism determined by $T_{\phi\psi}$,
$V_{\phi\psi}$ and $\widetilde{\phi\psi}$.  The following Lemma assures us that $\phi\psi$
is actually in $\HC (T,V)$.

\begin{lem}\label{lem:comp} The sets $T_{\phi\psi}$ and $V_{\phi\psi}$ defined above
are normal closed subsets
of $T$ and $V$.  The function $\widetilde{\phi\psi}_{X/T_{\phi\psi}}$ defined above is well-defined
and determines a based isomorphism of schemes $\widetilde{\phi\psi}:T\dm T_{\phi\psi}\to V\dm V_{\phi\psi}$.
Finally, we have
$t^{T_{\phi\psi}}\widetilde{\phi\psi}=v^{V_{\phi\psi}}$
if and only if there exists a $u\in U$ such that
$t^{T_{\phi}}\tilde\phi = u^{U_\phi}$ and $u^{U_\psi}\tilde \psi=v^{V_\psi}.$
\end{lem}

\begin{proof} We first recall two facts about quotients of schemes (see Lemmas 4.1.7 and 4.2.4
in Zieschang's work \cite{PHZ}).
First, whenever $R$ is a nonempty
subset of a scheme $S$ with a closed subset $T$, the quotient $R\dm T$ is closed if and only if
$R$ is closed.  Also, if $T$ and $U$ are closed subsets of $S$ with $T\subseteq U$,
then $U\dm T$ is normal in $S\dm T$ if and only if $U$ is normal in $S$.

Since $U_\phi$ and $U_\psi$ are normal and closed, then $U_\phi U_\psi$
is normal and closed, whence by the above comments, $(U_\phi U_\psi)\dm U_{\phi}$ is normal
and closed in $U\dm U_\phi$.
Now, observe that $T_{\phi\psi}$ contains $T_\phi$, and
$T_{\phi\psi}\dm T_\phi$ coincides precisely with the preimage of $(U_\phi U_\psi\dm U_{\phi})$
under $\tilde\phi$.  Since $\tilde \phi$ is an isomorphism, it follows that $T_{\phi\psi}\dm T_\phi$
is normal and closed in $T\dm T_\phi$, so $T_{\phi\psi}$ is normal and closed in $T$.
A similar argument shows that $V_{\phi\psi}$ is normal and closed.

Now, $T_{\phi\psi}$ is the kernel of the composition
\[\xymatrix{{T}\ar[r] & {T\dm T_\phi}\ar[r]^-{\tilde\phi} & {U\dm U_\phi}\ar[r] & {U\dm U_\phi U_\psi}}\]
where the first and third morphisms are the natural quotient morphisms.
There is then an induced isomorphism of schemes \[\bar\phi:T\dm T_{\phi\psi}\to U\dm U_\phi U_\psi.\]
The function $\bar\phi_{X/T_{\phi\psi}}$ takes an element
$xT_{\phi\psi}\in X/T_{\phi\psi}$ to the element $y (U_\phi U_\psi)\in Y/(U_\phi U_\psi)$, where
$y\in Y$ is any coset representative of
$(x T_\phi)\tilde \phi\in Y/ U_\phi$.
Likewise, the function $\bar\phi_{T\dm T_{\phi\psi}}$ takes
$t^{T_{\phi\psi}}\in T\dm T_{\phi\psi}$ to the element
$u^{U_\phi U_\psi}$ where $u^{U_{\phi}}=t^{T_\phi}\tilde\phi$.
Since $\tilde\phi$ is a based morphism, $\bar\phi$ preserves the basepoint.

Similarly, $V_{\phi\psi}$ is the kernel of the composition
\[\xymatrix{{V}\ar[r] & {V\dm V_\psi}\ar[r]^-{\tilde\psi^{-1}} & {U\dm U_\psi}\ar[r] & {U\dm U_\phi U_\psi}}.\]
There is again an induced isomorphism
\[\overline{\psi^{-1}}:V\dm V_{\phi\psi}\to U\dm U_\phi U_\psi.\]
If we then let $\bar\psi$ denote the inverse of $\overline{\psi^{-1}}$, then $\bar\psi_{Y/U_\phi U_\psi}$
takes an element $y (U_\phi U_\psi)$ to $z V_{\phi\psi}$, where
$z$ is any coset representative of $(y U_\psi)\tilde\psi$.
Likewise, $\bar\psi$ takes $u^{U_\phi U_\psi}$ to $v^{V_{\phi\psi}}$
where $u^{U_\psi}\tilde\psi=v^{V_\psi}$.
As above, $\bar\psi$ preserves the basepoint.

Now notice that
$(\bar\phi\bar\psi)_{X/T_{\phi\psi}}$ coincides with $\widetilde{\phi\psi}_{X/T_{\phi\psi}}$.
The second statement of the 
lemma therefore follows since $\bar\phi\bar\psi$ is a based isomorphism of schemes.
For the final statement, we have
$t^{T_{\phi\psi}}\widetilde{\phi\psi}=v^{V_{\phi\psi}}$
if and only if $t^{T_{\phi\psi}}\bar\phi\bar\psi=v^{V_{\phi\psi}}$.
If there is a $u\in U$ such that $t^{T_\phi}\tilde \phi=u^{U_\phi}$ and
$u^{U_\psi}\tilde\psi=v^{V_\psi}$, then by our above descriptions
of $\bar\phi_{T\dm T_{\phi\psi}}$ and
$\bar\psi_{U\dm U_\phi U_\psi}$, we have
$t^{T_{\phi\psi}}\bar\phi=u^{U_\phi U_\psi}$ and
$u^{U_\phi U_\psi}\bar\psi=v^{V_{\phi\psi}},$ so
indeed $t^{T_{\phi\psi}}\bar\phi\bar\psi=v^{V_{\phi\psi}}$.
On the other hand, if $t^{T_{\phi\psi}}\bar\phi\bar\psi=v^{V_{\phi\psi}}$,
then for some $u_0\in U$ and $v_0\in V$, we have $t^{T_{\phi}}\tilde\phi=u_0^{U_\phi}$
and $u_0^{U_\psi}\tilde\psi=v_0^{V_{\psi}}$, where $v_0^{V_{\phi\psi}}=v^{V_{\phi\psi}}$.
Since $V_{\phi\psi}$ is normal, then $v\in v_0v_1$ for some $v_1\in V_{\phi\psi}$.
By definition of $V_{\phi\psi}$, there is a $u_1\in U_\phi$ such that
$u_1^{U_\psi}\tilde\psi=v_1^{V_\psi}$.  Thus, $(u_0^{U_\psi}u_1^{U_\psi})\tilde\psi$
contains $v^{V_\psi}$, so we can find $u\in u_0u_1$ such that $u^{U_\psi}\tilde\psi=v^{V_\psi}$.
Since $u_1\in U_\phi$, we have $u^{U_\phi}=u_0^{U_\phi}$, so $t^{T_{\phi}}\tilde\phi=u^{U_{\phi}}$
and $u^{U_\psi}\tilde\psi=v^{V_{\psi}}$, as needed.
\end{proof}

The following lemma now implies that $\mathcal C$ as defined above is a category.

\begin{lem} The composition defined above is associative.  Also, if $T$ is a scheme on $X$,
and we let $T_{\id_T}$ be $\{1_X\}$ and $\tilde\id_T$ denote the identity
morphism on $T\dm \{1_X\}$, then $\id_T$ is the identity morphism for $T$ in $\Cc$.
\end{lem}

\begin{proof}  The second statement is easy, so we only prove the first.
We suppose that $S$, $T$, $U$ and $V$ are schemes on based sets $W$, $X$, $Y$, and $Z$
respectively, and $\phi$, $\chi$ and $\psi$ are morphisms in $\HC (S,T)$, $\HC(T,U)$, and
$\HC(U,V)$ respectively.  We must show that $(\phi\chi)\psi=\phi(\chi\psi)$.
Using the third statement of Lemma \ref{lem:comp}, it is straightforward to show that
both $S_{(\phi\chi)\psi}$ and $S_{\phi(\chi\psi)}$ coincide with the set
\[\{s\in S:\text{ for some }t\in T\text{ and }u\in U\text{ we have }
s^{S_\phi}\tilde\phi=t^{T_\phi},
t^{T_\chi}\tilde \chi=u^{U_\chi}\text{ and }
u^{U_\psi}\tilde \psi=1_Z^{V_\psi}\}.\]
Thus, $S_{(\phi\chi)\psi}=S_{\phi(\chi\psi)}$.  By a similar argument,
$V_{(\phi\chi)\psi}=V_{\phi(\chi\psi)}$.  Finally, we must show
$\widetilde{(\phi\chi)\psi}=\widetilde{\phi(\chi\psi)}$, and for this, it
suffices to show that $\widetilde{(\phi\chi)\psi}_{W/S_{(\phi\chi)\psi}}=
\widetilde{\phi(\chi\psi)}_{W/S_{\phi(\chi\psi)}}$.  But given $w\in W$,
we may choose $x\in X$, $y\in Y$, and $z\in Z$ such that
$(wS_\phi)\tilde \phi=xT_\phi$,
$(xT_\chi)\tilde \chi=yU_\chi$, and
$(yU_\psi)\tilde \psi=zV_\psi$.
Then by the definition we gave immediately before the statement of Lemma \ref{lem:comp}, both $\widetilde{(\phi\chi)\psi}_{W/S_{(\phi\chi)\psi}}$ and
$\widetilde{\phi(\chi\psi)}_{W/S_{\phi(\chi\psi)}}$ take
$w$ to $zV_{(\phi\chi)\psi}=zV_{\phi(\chi\psi)}$. 
\end{proof}

Given two based schemes $T$ and $U$ on based sets $X$ and $Y$,
the set $\HC(T,U)$ can be given the structure of a partially ordered set.
Given two morphisms $\phi,\psi\in \HC(T,U)$, we will say that $\phi\leq \psi$ if
$T_\phi\subseteq T_\psi$, $U_\phi\subseteq U_\psi$, and for each
$x\in X$, $(xT_\phi)\tilde\phi\subseteq (xT_\psi)\tilde \psi$.

\begin{notn}\label{notn:star} If $\phi\in \HC(T,U)$, then we have a
morphism $\phi^*\in \HC(U,T)$ defined as follows:
$T_{\phi^*}=T_{\phi}$, $U_{\phi^*}=U_{\phi}$, and $\widetilde{\phi^*}=\tilde{\phi}^{-1}$.
\end{notn}

Note that if $\phi\in \HC(T,U)$, then the composition $\phi\phi^*\in \HC(T,T)$ will
not typically be the identity.  Instead, both normal closed subsets of $\phi\phi^*$
coincide with $T_\phi$, and $\widetilde{\phi\phi^*}$ is the identity morphism
on the scheme $T\dm T_\phi$.  Thus, $\id_T\leq \phi\phi^*$.

\section{Labelling sets}\label{sec:label}

\begin{defn}\label{def:wlabel} A {\it weak labelling set} is a set $\tau$ equipped with
an involution $*:\tau\to \tau$, a distinguished element $1\in \tau$ fixed by $*$,
and for each triple $p, q, r\in \tau$, a nonnegative integer
$a^\tau_{pqr}$.  We write $p^*$ for the image of $p$ under $*$. 
\end{defn}

\begin{defn}\label{def:tauscheme}
Given a weak labelling set $\tau$, a {\it $\tau$-scheme} is a pair $(T,\al)$, where
$T$ is a scheme on a based set $X$, and $\al:\tau\to T$ is a bijection satisfying the following:
\end{defn}

\begin{enumerate}
	\item $1\al=1_X$,
	\item for each $p\in \tau$, $(p^*)\al=(p\al)^*$,
	\item for each triple $p, q, r\in \tau$, $a^\tau_{pqr}=a_{(p\al)(q\al)(r\al)}$.
\end{enumerate}	

\begin{defn}\label{def:label}
A {\it labelling set} is a weak labelling set $\tau$ for which there exists a $\tau$-scheme.
\end{defn}

For example, a scheme $T$ on a set $X$ determines a labelling set $\tau_T$
by forgetting about the underlying set $X$.  That is, the set $\tau_T$ is equal to the
set $T$, the involution in $\tau_T$ is equal to the involution on $T$, and the distinguished
element $1\in \tau_T$ is equal to $1_X$.  Finally, if $p, q, r\in \tau_T=T$, then
$a^\tau_{pqr}$ is the corresponding structure constant for $T$.
As another example, if $S$ is a scheme on a set $X$, and $T$ is a closed subset of $S$,
and $\tau=\tau_T$, then for each $x\in X$, the subscheme of $S$ defined
by $xT$ is a $\tau$-scheme:  the bijection $\al$ takes $t\in \tau_T=T$ to $t\cap (xT\times xT)$.

\begin{rem}\label{rem:closure}
If $(T,\al)$ and $(U,\beta)$ are both $\tau$-schemes, and $r, p, q\in \tau$, then
$r\al\in (p\al)(q\al)$ if and only if $r\beta\in (p\beta)(q\beta)$.  Indeed,
$r\al\in (p\al)(q\al)$ if and only if $a_{(p\al)(q\al)(r\al)}>0$,
which holds if and only if $a^\tau_{pqr}>0$.  Similarly, $r\beta\in (p\beta)(q\beta)$
if and only if $a^\tau_{pqr}>0$.
\end{rem}

\begin{defn}\label{def:mult}
Given a labelling set $\tau$, and given $P\subseteq \tau$ and $Q\subseteq \tau$, we let
\[PQ=\{r\in\tau:a^\tau_{pqr}>0\text{ for some }p\in P, q\in Q\}.\]
If $\tau'$ is a subset of $\tau$ and $p,q\in \tau$, we define $p\tau'=\{p\}\tau'$,
$\tau'q=\tau'\{q\}$ and $pq=\{p\}\{q\}$.
\end{defn}

To see that the product on subsets of $\tau$ is associative,
choose some $\tau$-scheme $(T,\al)$.  Then by Condition (3) of Definition
\ref{def:tauscheme}, we have $(PQ)\al=(P\al)(Q\al)$.  Since the complex product for association schemes
is associative, we have for any $P, Q, R\subseteq \tau$
\[((PQ)R)\alpha=((P\alpha)(Q\alpha))R\alpha=P\alpha((Q\alpha)(R\alpha))=(P(QR))\alpha,\]
so $(PQ)R=P(QR)$ since $\alpha$ is a bijection.

\begin{defn}\label{def:ncl}
A subset $\tau'\subseteq \tau$ is {\it closed} if $pq^*\subseteq \tau'$
whenever $p,q\in \tau'$.  A subset $\tau'\subseteq \tau$ is
{\it normal} if $p\tau'=\tau' p$ whenever $p\in \tau$.
\end{defn}

\begin{rem}\label{rem:ncl} If $(T,\al)$ is any $\tau$-scheme, then $\alpha$ induces a one-to-one
correspondence between the closed subsets
of $\tau$ and the closed subsets of $T$.  Indeed, if $\tau'$ is closed, and
$p\al, q\al$ are arbitary elements in $\tau'\al$, then
as in Remark \ref{rem:closure}, for any $r\al\in (p\al) (q\al)^*=(p\al) (q^*\al)$,
we must have $a^\tau_{pq^*r}>0$, so $r\in p^*q$.  This implies $r\in \tau'$, so $r\al\in \tau'\al$,
and thus $\tau'\al$ is closed.  Conversely, if $\tau'\al$ is closed, and $p,q\in \tau'$, then
for any $r\in pq^*$, we have $r\alpha\in
(p\alpha)(q^*\alpha)=(p\alpha)(q\alpha)^*\subseteq \tau'\al$, so $r\in \tau'$, whence $\tau'$ is closed.

By a similar argument,
if $(T,\al)$ is a $\tau$-scheme, then $\alpha$ induces a one-to-one correspondence between
the normal subsets of $\tau$ and the normal subsets of $T$.
Thus, if $\tau_1$ and $\tau_2$ are two normal closed subsets of a labelling set $\tau$, then
$\tau_1\tau_2$ must be a normal closed subset of $\tau$.  To see this, just choose a $\tau$-scheme
$(T,\al)$; then the image of $\tau_1\tau_2$ under $\al$ corresponds to $(\tau_1\al)(\tau_2\al)$,
which is a normal closed subset of $T$.

Finally, if $\tau'$ is closed and normal in $\tau$, and $(T,\al)$ is a $\tau$-scheme, then 
the partition of $T$ induced by $\tau'\al$ corresponds under $\al$ to a partition of $\tau$.
Thus, the cosets $p\tau'$, with $p\in\tau$, form a partition of $\tau$.
\end{rem}

\begin{notn} Let $\Cc(\tau)$ denote the category whose objects are $\tau$-schemes, and whose
morphisms are morphisms in the underlying category $\Cc$.  That is, if $(T,\al)$ and $(U,\beta)$
are $\tau$-schemes, then $\HCT((T,\al),(U,\beta))=\HC(T,U)$.
\end{notn}

In the following definition, we let $\Pp(\tau)$ denote the power set of $\tau$.

\begin{defn}\label{def:phitau}
Given $\phi\in \HCT((T,\al),(U,\beta))$, let $\phi(\tau):\Pp(\tau)\to \Pp(\tau)$ be the function
taking $P\subseteq \tau$ to
\[\{u\in \tau:\text{ for some }t\in P, (t\al)^{T_\phi} \tilde\phi=(u\beta)^{U_\phi}.\}\]
\end{defn}

\begin{rem}\label{rem:recover}
In general, given a function $\al:\Pp(\tau)\to \Pp(\tau)$, we define
\[\ker(\al)=\{t\in \tau:\{t\}\al=\{1\}\al\}.\]
Then if $\phi\in \HCT((T,\al),(U,\beta))$, it is easy to check that $U_\phi=\left(\{1\}\phi(\tau)\right)\beta$ and
$T_\phi=\ker(\phi(\tau))\al.$
That is, we can recover $T_\phi$ and $U_\phi$ from $\phi(\tau)$,
$\al$, and $\beta$.
\end{rem}

\section{Semidirect products from actions}\label{sec:act}

\begin{defn}\label{defn:action} Suppose $T$ is a scheme on a based set $X$, $U$ a scheme
on a based set $Y$, and $\tau=\tau_T$ is the
labelling set corresponding to $T$.
Then an {\it action $\zeta$ of $U$ on $T$} consists of
a $\tau$-scheme $\zeta_y=(T_y,\al^y)$ on $X$
for each $y\in Y$, and a morphism $\zeta_{y_1}^{y_2}\in \HC(T_{y_1},T_{y_2})$ for each pair $y_1, y_2\in Y$.
We require the following properties to hold

\begin{enumerate}
	\item If $y_*\in Y$ is the basepoint in $Y$, then $T_{y_*}=T$ and $\al^{y_*}:\tau_T\to T_{y_*}=T$
	is the identity.
	(Recall that the underlying set of $\tau_T$ is defined to be $T$.) 
	\item For each $y\in Y$, $\zeta_y^y$ is the identity morphism on $T_y$.
	\item For each pair $y_1, y_2\in Y$, $\zeta_{y_2}^{y_1}=(\zeta_{y_1}^{y_2})^*$
	(see Notation \ref{notn:star}).
	\item The function $\zeta_{y_1}^{y_2}(\tau):\Pp(\tau)\to \Pp(\tau)$ depends
	only on the scheme element $u\in U$ containing $(y_1,y_2)$.
	\item For each triple $y_1, y_2, y_3\in Y$, we have
	$\zeta_{y_1}^{y_3}\leq \zeta_{y_1}^{y_2}\zeta_{y_2}^{y_3}.$
\end{enumerate}
\end{defn}

\begin{notn}\label{notn:simplify}
We will often denote $(T_{y_1})_{\zeta_{y_1}^{y_2}}$
as $T'_{y_1y_2}$ and $(T_{y_2})_{\zeta_{y_1}^{y_2}}$ as $T''_{y_1y_2}$.
Thus, $\tilde\zeta_{y_1}^{y_2}$
is an isomorphism of schemes from $T_{y_1}\dm T'_{y_1y_2}$ to $T_{y_2}\dm T''_{y_1y_2}$.
We will let $\zeta_u=\zeta_{y_1}^{y_2}(\tau)$, for some $(y_1, y_2)\in u$.
By condition (4), $\zeta_u$ is well-defined.  Let $\tau'_u=\ker \zeta_u$
and $\tau''_u=\{1\}\zeta_u$.  Note that by Remark \ref{rem:recover},
the normal closed subsets $T'_{y_1y_2}$ and $T''_{y_1y_2}$
coincide with $(\tau'_u)\al^{y_1}$ and $(\tau''_u)\al^{y_2}$.
\end{notn}

Given an action $\zeta$ of $U$ on $T$, where $U$ is a based scheme on $Y$
and $T$ is a based scheme on $X$, we will define a corresponding scheme 
$U\ltimes_\zeta T$ on the set $Y\times X$.  
Given $u\in U$ and $t\in \tau=\tau_T$, we let
$[u,t]$ denote the set of all pairs $((y_1, x_1), (y_2, x_2))\in (Y\times X)^{\times 2}$
such that $(y_1,y_2)\in u$ and \[((x_1 T'_{y_1y_2})
\tilde\zeta_{y_1}^{y_2}, x_2 T''_{y_1y_2})
\in (t\al^{y_2})^{T''_{y_1y_2}}\]

\begin{lem}\label{lem:partition} As $u$ ranges over $U$ and $t$ over $\tau$,
the sets $[u,t]$ form a partition of $(Y\times X)^{\times 2}$.
\end{lem}

\begin{proof}
Given an arbitrary pair $((y_1,x_1),(y_2,x_2))$, we have
$(y_1,y_2)\in u$ for some $u\in U$, and
 \[((x_1 T'_{y_1y_2})
\tilde\zeta_{y_1}^{y_2}, x_2 T''_{y_1y_2})
\in (t\al^{y_2})^{T''_{y_1y_2}}\] 
for some $t\in \tau$, since $\al^{y_2}:\tau\to T_{y_2}$
is a bijection.  Thus every element in $(X\times Y)^{\times 2}$ belongs
to $[u,t]$ for some $u\in U, t\in \tau$.

If $((y_1,x_1),(y_2,x_2))\in [u,t]\cap [\hat u,\hat t]$, then
$(y_1,y_2)\in u\cap \hat u$, so $u=\hat u$ since $U$ is a partition
of $Y\times Y$.  Similarly
$(t\al^{y_2})^{T''_{y_2}}=(\hat t\al^{y_2})^{T''_{y_2}}$; equivalently
(since $T''_{y_2}$ is normal and $\al^{y_2}$ is a bijection),
we can find $s\in T''_{y_1y_2}(\al^{y_2})^{-1}$ such that
$t\al^{y_2}\in (\hat t\al^{y_2})(s\al^{y_2}).$
We claim $[u,t]=[u,\hat t]$.
If $((y_3,x_3),(y_4,x_4))\in [u,t]$, then $(y_3,y_4)\in u$.  Also,
since \[s\in T''_{y_1y_2}(\al^{y_2})^{-1}=\tau''_u=T''_{y_3y_4}(\al^{y_4})^{-1},\]
we have 
$s\al^{y_4}\in T''_{y_3y_4}$.  By Remark \ref{rem:closure},
$t\al^{y_2}\in (\hat t\al^{y_2})(s\al^{y_2})$ implies 
that $t\al^{y_4}\in (\hat t\al^{y_4})(s\al^{y_4}),$ so
$(t\al^{y_4})^{T''_{y_3y_4}}=(\hat t\al^{y_4})^{T''_{y_3y_4}}$, which in
turn implies that $((y_3,x_3),(y_4,x_4))\in [u,\hat t]$.  Thus,
$[u,t]\subseteq [u,\hat t]$, and by a symmetric argument, $[u,\hat t]\subseteq [u,t]$.
\end{proof}

\begin{notn} If $\zeta$ is an action of $U$ on $T$, where $U$ and $T$
are schemes on based sets $Y$ and $X$, we let $U\ltimes_\zeta T$ denote the set of all
$[u,t]\subseteq (Y\times X)^{\times 2}$
as $u$ ranges over $U$ and $t$ over $\tau=\tau_T$.
\end{notn}

\begin{lem}\label{lem:1} The set $U\ltimes_\zeta T$ contains
\[1_{Y\times X}=\{((y,x),(y,x)):y\in Y\text{ and }x\in X\}.\]
\end{lem}

\begin{proof}
We claim $[1_Y,1]=1_{Y\times X}$, where $1$ denotes the distinguished element in $\tau$.
Given any $((y_1,x_1),(y_2,x_2))\in (Y\times X)^{\times 2}$, we have $(y_1,y_2)\in 1_Y$ if and only if
$y_1=y_2$.  By Condition (1) in Definition
\ref{defn:action}, $\zeta_{y_1}^{y_1}$ is the identity morphism in $\HC(T_{y_1},T_{y_1})$.
Thus, $T'_{y_1y_1}=T''_{y_1y_1}=1_{\{X\}}$, and $\tilde \zeta_{y_1}^{y_1}$ is the identity morphism
of schemes.  Since $1\al^{y_2}=1_X$ by Condition (1) of Definition \ref{def:wlabel},
$((y_1,x_1),(y_2,x_2))\in [1_Y,1]$ if and only if $y_1=y_2$ and $x_1=x_2$.
\end{proof}

\begin{lem}\label{lem:*} Suppose $[u,t]\in U\ltimes_\zeta T$.  Let $\hat t^*$ be any element in 
$\{t^*\}\zeta_{u^*}\subseteq \tau$.  Then
$[u,t]^*=[u^*,\hat t^*]$.  In particular, $[u,t]^*\in U\ltimes_\zeta T$.
\end{lem}

\begin{proof}
First, $((y_1, x_1), (y_2,x_2))\in [u,t]$ if and only if
$(y_1,y_2)\in u$ and \[((x_1 T'_{y_1y_2})
\tilde\zeta_{y_1}^{y_2}, x_2 T''_{y_1y_2})
\in (t\al^{y_2})^{T''_{y_1y_2}}.\]
Equivalently, $(y_2,y_1)\in u^*$ and (by Condition 2 of Definition \ref{def:wlabel}),
\[(x_2 T''_{y_1y_2},
(x_1 T'_{y_1y_2})
\tilde\zeta_{y_1}^{y_2})
\in (t^*\al^{y_2})^{T''_{y_1y_2}}.\]
Since $\tilde\zeta_{y_1}^{y_2}$ is an isomorphism of schemes, this is equivalent to saying
$(y_2,y_1)\in u^*$ and
\[((x_2 T''_{y_1y_2})(\tilde\zeta_{y_1}^{y_2})^{-1},
x_1 T'_{y_1y_2})
\in \left((t^*\al^{y_2})^{T''_{y_1y_2}}\right)(\tilde\zeta_{y_1}^{y_2})^{-1}.\]
By Condition (3) of Definition \ref{defn:action},
$\zeta_{y_2}^{y_1}=(\zeta_{y_1}^{y_2})^*$, so the above is equivalent to saying $(y_2,y_1)\in u^*$ and
\[((x_2 T'_{y_2y_1})\tilde\zeta_{y_2}^{y_1},
x_1 T''_{y_2y_1})
\in \left((t^*\al^{y_2})^{T'_{y_2y_1}}\right)(\tilde\zeta_{y_2}^{y_1}).\]
From the definition of $\hat t^*$ and $\zeta_{u^*}$, we have
\[\left((t^*\al^{y_2})^{T'_{y_2y_1}}\right)(\tilde\zeta_{y_2}^{y_1})=
(\hat t^*\al^{y_1})^{T''_{y_2y_1}}.\]
Thus, $((y_1,x_1),(y_2,x_2))\in [u,t]$ if and only if $((y_2,x_2),(y_1,x_1))\in [u^*,\hat t^*]$.
\end{proof}

In order to prove the third scheme condition for $U\ltimes_\zeta T$, 
we need the following preliminary lemma.

\begin{lem}\label{lem:prelim1}
Suppose that $\tau$ is a labelling set, $\tau_1$ and $\tau_2$ are normal closed subsets of $\tau$
(see Definition \ref{def:ncl}), and $p, q, r$ are elements in $\tau$.  If $(T,\al)$ is any $\tau$-scheme
and $(y,z)\in r\al$, then $|y(p\al)(\tau_1\al)\cap z(q\al)^*(\tau_2\al)|$ depends
only on $\tau_1$, $\tau_2$, $p$, $q$ and $r$.  That is, this number does not
depend on $(T, \al)$ or $(y,z)$.
\end{lem}

\begin{proof}
First, let $T_1=\tau_1\al$ and $T_2=\tau_2\al$.  As observed in Remark \ref{rem:ncl}, $T_1$ and
$T_2$ are normal closed subsets of $T$.  We have $x\in y(p\al)T_1\cap z(q\al)^*T_2$ 
if and only if $x\in yp'\alpha\cap z(q'\alpha)^*$ for some $p'\al\in (p\al)T_1$
and some $q'\al\in ((q\al)^*T_2)^*$.
Since $T_2$ is closed and normal, $((q\al)^*T_2)^*=T_2(q\al)=(q\al)T_2$.
Therefore, \[|y(p\al)T_1\cap z(q\al)^*T_2|=
\sum_{p'\al\in (p\al)T_1}\sum_{q'\al\in (q\al)T_2} a_{(p'\al) (q'\al) (r\al)}
=\sum_{p'\in p\tau_1}\sum_{q'\in q\tau_2} a^\tau_{p'q'r}.\]
This number depends only on $\tau_1$, $\tau_2$, $p$, $q$ and $r$.
\end{proof}

\begin{lem}\label{lem:apqr} Suppose $p_i=[u_i,t_i]\in U\ltimes_\zeta T$ for $i=1, 2, 3$.  Then
there is a nonnegative integer $a_{p_1p_2p_3}$ such that for any $((y_1,x_1),(y_2,x_2))\in p_3$,
$|(y_1,x_1)p_1\cap (y_2,x_2)p_2^*|=a_{p_1p_2p_3}.$
\end{lem}

\begin{proof}
First, suppose $u_3\notin u_1u_2$.  Then given $((y_1,x_1),(y_2,x_2))\in p_3$, 
$(y_1,x_1)p_1\cap (y_2,x_2)p_2^*$ must be empty:  if it contained $(y,x)$, then
we would have $(y_1,y_2)\in u_3$, $(y_1,y)\in u_1$ and $(y,y_2)\in u_2$ so $u_3\in u_1u_2$.
Thus, if $u_3\notin u_1u_2$, then $a_{p_1p_2p_3}=0$.
We may therefore just consider the case when $u_3\in u_1u_2$.  

Suppose given $((y_1,x_1),(y_2,x_2))\in p_3$.  Fix an element $y\in y_1u_1\cap y_2u_2^*$.
We will show that the cardinality
\begin{equation}
\label{eqn:set}
\left|\{x\in X:(y,x)\in (y_1,x_1)p_1\cap (y_2,x_2)p_2^*\}\right|
\end{equation}
depends only on $p_1, p_2,$ and $p_3$.
Since this cardinality is independent of $y$, multiplying it by $a_{u_1u_2u_3}=|y_1u_1\cap y_2u_2^*|$
then yields the needed value for $a_{p_1p_2p_3}$.

Since $((y_1, x_1), (y_2, x_2))\in [u_3,t_3]$, we have
\[\left((x_1T'_{y_1y_2})\tilde \zeta_{y_1}^{y_2},
x_2T''_{y_1y_2}\right)\in (t_3\al^{y_2})^{T''_{y_1y_2}}.\]
Since $T''_{y_1y_2}$ is normal, we may choose a coset representative $\bar x_1$ of
$(x_1T'_{y_1y_2})\tilde \zeta_{y_1}^{y_2}$ such that
$(\bar x_1, x_2)\in t_3\al^{y_2}$, so \begin{equation}\label{eqn:A}
(\bar x_1 T'_{y_2y},x_2 T'_{y_2y})\in (t_3\al^{y_2})^{T'_{y_2y}}.\end{equation}
Since $\zeta_{y_1}^{y_2}\leq \zeta_{y_1}^y\zeta_{y}^{y_2}$, the scheme isomorphism
$\widetilde{\zeta_{y_1}^{y}\zeta_{y}^{y_2}}$ takes the coset of $x_1$ to the coset of $\bar x_1$,
so by the definition of $\widetilde{\zeta_{y_1}^{y}\zeta_{y}^{y_2}}$ given in the paragraph preceding
Lemma \ref{lem:comp}, we can choose $\hat x_1$ such that
\begin{equation}\label{eqn:B}
(x_1T'_{y_1y})(\tilde \zeta_{y_1}^y)=\hat x_1 T''_{y_1y}
\end{equation}
and
\begin{equation}\label{eqn:C}
(\hat x_1 T'_{yy_2})(\tilde \zeta_y^{y_2})=\bar x_1 T''_{yy_2}.
\end{equation}
Since $\zeta_{y}^{y_2}=(\zeta_{y_2}^y)^*$, we have $T''_{yy_2}=T'_{y_2y}$, and
\begin{equation}\label{eqn:D}
\hat x_1 T'_{yy_2}=(\bar x_1 T'_{y_2y})(\tilde \zeta_{y_2}^y).
\end{equation}
Since $(y_2,y)\in u_2^*$,
it follows from the definition of $\zeta_{u_2}^*$ that for any
$r\in \{t_3\}\zeta_{u_2^*}$, we have
\[(r\alpha^y)^{T''_{y_2y}}=(t_3\alpha^{y_2})^{T'_{y_2y}}\tilde\zeta_{y_2}^y,\]
which, together with Equation \ref{eqn:A} yields
\[((\bar x_1T'_{y_2y})(\tilde \zeta_{y_2}^y),(x_2T'_{y_2y})(\tilde \zeta_{y_2}^y))
\in (r\al^y)^{T''_{y_2y}}=(r\al^y)^{T'_{yy_2}}.\]
Thus from Equation \ref{eqn:D},
\[(\hat x_1 T'_{yy_2},(x_2T'_{y_2y})(\tilde \zeta_{y_2}^y))
\in (r\al^y)^{T'_{yy_2}}.\]
Since $T'_{yy_2}$ is normal, we may choose a coset representative
$\hat x_2\in  (x_2T'_{y_2y})(\tilde \zeta_{y_2}^y)$
such that $(\hat x_1,\hat x_2)\in r\al^y$.

Now, recalling that $y$ is fixed, we have $(y,x)\in (y_1,x_1)p_1$ if and only if
\[((x_1T'_{y_1y})\tilde\zeta_{y_1}^y,xT''_{y_1y_2})\in (t_1\alpha^y)^{T''_{y_1y}}.\]
From Equation \ref{eqn:B}, this holds if and only if
\[(\hat x_1T''_{y_1y},xT''_{y_1y})\in (t_1\alpha^y)^{T''_{y_1y}}.\]
Since $T''_{y_1y}$ is normal in $T_y$, this holds if and only if
$x\in \hat x_1(t_1\al^y)T''_{y_1y}.$
By a similar argument, together with Lemma \ref{lem:*} and Condition (2) of \ref{def:tauscheme},
we have $(y,x)\in (y_2, x_2)p_2^*$ if and only if
$x\in \hat x_2(\hat t_2\al^y)^*T''_{y_2y}$, where
$\hat t_2^*\in \tau$ belongs to $\{t_2^*\}\zeta_{u_2^*}$.
Now, recall from Notation \ref{notn:simplify}
that $T''_{y_1y}=\tau''_{u_1}\al^y$ and $T''_{y_2y}=T'_{yy_2}=\tau'_{u_2}\al^{y}$.
By Lemma \ref{lem:prelim1},
\begin{equation}\label{eqn:set2}
|\hat x_1(t_1\al^y)(\tau''_{u_1}\al^y)\cap \hat x_2(\hat t_2\al^y)^*(\tau'_{u_2}\al^{y})|
\end{equation}
depends only on $t_1, \hat t_2, \tau''_{u_1}, \tau'_{u_2}$ and $r$, hence only on 
$t_1, \hat t_2, u_1, u_2$ and $r$.  But since $\hat t_2$ is any element in
$(\{t_2^*\}\zeta_{u_2^*})^*$
and $r$ is any element in
$\{t_3\}\zeta_{u_2^*}$, $\hat t_2$ and $r$ depend only on $t_2$, $t_3$ and $u_2$.
Thus, the cardinality of the set (\ref{eqn:set2}), and hence of set (\ref{eqn:set}),
depends only on $p_1$, $p_2$ and $p_3$.
\end{proof}

Now, if $Y$ and $X$ are based sets, then $Y\times X$ has a natural basepoint $(y_*,x_*)$.
The following corollary now follows from Lemmas \ref{lem:1}, \ref{lem:*}, and \ref{lem:apqr}.

\begin{cor} The set $U\ltimes_\zeta T$ is a based scheme on $Y\times X$.
\end{cor}

\section{Properties of the semidirect product}\label{sec:prop}

Suppose $U$ and $T$ are based schemes on $X$ and $Y$, $\tau=\tau_T$, and
$\zeta$ is an action of $U$ on $T$.  In this section, we show that the scheme
$U\ltimes_\zeta T$ contains a normal closed subset $\tilde T:=\{[1_Y,t]:t\in \tau\}$, such that
the subscheme of $U\ltimes_\zeta T$ defined by $(y_*,x_*)\tilde T$ is isomorphic
to $T$.  We also produce a morphism of schemes $i:U\to U\ltimes_\zeta T$ which splits
the natural quotient morphism $\pi:U\ltimes_\zeta T\to (U\ltimes_\zeta T)\dm \tilde T$.
This justifies calling $U\ltimes_\zeta T$ a semidirect product.  We also show that
Condition \ref{con:hyp} holds.

\begin{prop}\label{prop:cn} The subset $\tilde T$ is a normal closed subset of $U\ltimes_\zeta T$.
Moreover, the subscheme of $U\ltimes_\zeta T$ defined by $(y_*,x_*)\tilde T$ is isomorphic to $T$.
\end{prop}

\begin{proof} First, given $[1,t]\in \tilde T$, we have
$((y_1,x_1),(y_2,x_2))\in [1,t]$ if and only if $y_2=y_1$ and $(x_1,x_2)\in t\al^{y_1}$,
since $\zeta_{y_1}^{y_1}$ is the identity.  It follows from Definition \ref{def:tauscheme}
and Definition \ref{def:mult}
that $[1,t_1][1,t_2]=\cup_{t\in t_1t_2}[1,t]$ and $[1,t_1]^*=[1,t_1^*]$, so $\tilde T$ is closed.
Now, to show that $\tilde T$ is normal, it suffices to show
$\tilde T [u,t]\subseteq [u,t]\tilde T$ for any $[u,t]\in U\ltimes_\zeta T$.
Suppose given $(y_1,x_1), (y_2, x_2), (y_3, x_3)$, where $((y_1,x_1),(y_2,x_2))\in [1,t']$ for
some $t'\in \tau$, and $((y_2,x_2),(y_3,x_3))\in [u,t]$.  Then $y_1=y_2$ and $(y_2, y_3)\in u$,
so $(y_1, y_3)\in u$.  Also $(x_1,x_2)\in t'\al^{y_2}$ and
\[((x_2T'_{y_2y_3})\tilde \zeta_{y_2}^{y_3},
x_3T''_{y_2y_3})\in (t\al^{y_3})^{T''_{y_2y_3}}.\]
Now let $y_4=y_3$, and choose $x_4$ so that
\[((x_1T'_{y_1y_4})\tilde \zeta_{y_1}^{y_4},
x_4T''_{y_1y_4})\in (t\al^{y_4})^{T''_{y_1y_4}}.\]
Then $((y_1,x_1),(y_4,x_4))\in [u,t]$.  Choose $\tilde t\in \tau$ so that
$(x_4,x_3)\in \tilde t\al^{y_4}$.  Then \[((y_4,x_4),(y_3,x_3))\in [1,\tilde t]\in \tilde T.\]
Thus, $\tilde T[u,t]\subseteq [u,t]\tilde T$.

For the final statement, note that $(y_*,x_*)\tilde T=\{y_*\}\times X\subseteq Y\times X$.
Define $f_X:X\to \{y_*\}\times X$ by $xf_X=(y_*,x)$.
If $(x_1,x_2)\in t\in T$, then it follows
from Conditions (1) and (2) of Definition \ref{defn:action} that
\[(x_1f_X,x_2f_X)=((y_*,x_1),(y_*,x_2))\in [1,t\alpha^{y_*}]=[1,t].\]
Thus, $f_X$ determines a morphism of schemes $f:T\to \tilde T_{(y_*,x_*)\tilde T}$;
in particular, $tf_T=[1,t]$ for any $t\in T=\tau_T$.  Clearly, $f_X$ is a bijection, and $f_T$
is surjective.  If $tf_T=t'f_T$, then $[1,t]=[1,t']$, so for some $x_1, x_2\in X$, we have
$((y_*,x_1),(y_*,x_2))\in [1,t]\cap [1,t']$.  Thus, $(x_1,x_2)\in t\cap t'$, so $t=t'$, so $f_T$
is injective.
\end{proof}

Now define $i_Y:Y\to Y\times X$ by $yi_Y=(y,x_*)$, and let
$\pi:U\ltimes_\zeta T\to (U\ltimes_\zeta T)\dm \tilde T$ denote the natural quotient morphism.

\begin{prop}\label{prop:split}
The function $i_Y$ defines a based morphism of schemes $i:U\to U\ltimes_\zeta T$, and
$i \pi$ is an isomorphism of schemes.
\end{prop}

\begin{proof}  
Suppose $(y_1, y_2)\in u$.
Since $\tilde\zeta_{y_1}^{y_2}$ is based, we have
$(x_* T'_{y_1y_2})\tilde\zeta_{y_1}^{y_2}=x_* T''_{y_1y_2}$.  Since $1\al^{y_2}=1_X$,
we then have
\[((x_* T'_{y_1y_2})\tilde\zeta_{y_1}^{y_2},x_* T''_{y_1y_2})\in (1\al^{y_2})^{T''_{y_1y_2}}\]
Therefore, $((y_1,x_*),(y_2,x_*))\in [u,1]$.  Since $[u,1]$ does not depend on $(y_1, y_2)$,
and since $y_* i_Y=(y_*,x_*)$, $i_Y$ determines a based morphism $i:U\to U\ltimes_\zeta T$ of schemes.
Note that if $u\in U$, then $u i_U=[u,1]$.

Now, to show that $i \pi$ is an isomorphism, it suffices to show that $(i \pi)_Y$ is a bijection
and $(i \pi)_U$ is an injection, since if $(i \pi)_Y$ is surjective, then $(i\pi)_U$ is immediately surjective.
If $y_1 i \pi=y_2 i \pi$, then $(y_1,x_*)\tilde T=(y_2,x_*)\tilde T$, so $(y_1, x_*)\in (y_2, x_*)[1,t]$
for some $t\in \tau$, implying $y_1=y_2$.  On the other hand, given a coset $(y,x)\tilde T$, we have
$(y,x)\in (y,x_*)[1,t]$ for some $t\in \tau$, so $(y,x)\tilde T=(y,x_*)\tilde T=yi\pi$.
Thus, $(i\pi)_Y$ is a bijection.

Finally, suppose $u (i\pi)_U=v (i\pi)_U$, where $u,v\in U$.  Choose $(y_1, y_2)\in u$.  
Then \[((y_1, x_*)\tilde T,(y_2,x_*)\tilde T)\in [u,1]^{\tilde T}=[v,1]^{\tilde T}.\]
Thus, since $\tilde T$ is normal, there is some $t\in \tau$ so that
\[((y_1, x_*),(y_2, x_*)\in [v,1][1_Y,t].\]
That is, we can find some $(y,x)$ so that
\[((y_1, x_*), (y,x))\in [v, 1]\text{ and }((y,x),(y_2,x_*)\in [1_Y,t].\]
Therefore $(y_1, y_2)=(y_1,y)\in v$.
But since $(y_1, y_2)\in u$, we must have $u=v$, so $(i\pi)_U$ is injective.
\end{proof}

The following proposition shows that $U\ltimes_\zeta T$ satisfies
Condition \ref{con:hyp}.

\begin{prop} For each $u\in U$ and each $\tilde t\in \tilde T$, $|(ui)\tilde t|=|\tilde t(ui)|=1$.
\end{prop}

\begin{proof} We must show $|[u,1][1_Y,t]|=|[1_Y,t][u,1]|=1$ whenever $u\in U$ and $t\in \tau$.
We first show that $[u,1][1_Y,t]=[u,t]$.
Suppose $((y_1, x_1),(y_2,x_2))\in [u,1]$ and $((y_2,x_2),(y_3,x_3))\in [1_Y,t]$.
Then, we must have $y_3=y_2\in y_1u$ and (since $\zeta_{y_2}^{y_3}=\zeta_{y_2}^{y_2}$
is the identity morphism),
$x_3\in x_2(t\al^{y_3})$.  Also,
\[x_2T''_{y_1y_2}=(x_1T'_{y_1y_2})\tilde\zeta_{y_1}^{y_2}=(x_1T'_{y_1y_3})\tilde\zeta_{y_1}^{y_3}.\]
Thus,
\[x_3T''_{y_1y_3}=x_3T''_{y_1y_2}\in (x_2T''_{y_1y_2})(t\al^{y_3})^{T''_{y_1y_2}}=
(x_1T'_{y_1y_3})\tilde\zeta_{y_1}^{y_3}(t\al^{y_3})^{T''_{y_1y_3}}.\]
Thus, $((x_1, y_1),(x_3,y_3))\in [u,t]$, so $[u,1][1_Y,t]=[u,t]$.  That is,
$|(ui)\tilde t|=1$ for all $u\in U$ and $\tilde t\in \tilde T$.

Therefore, given $u\in U$ and $\tilde t\in \tilde T$, we have
$(\tilde t(ui))^*=(ui)^*\tilde t^*=(u^*i)\tilde t^*$.  By Proposition \ref{prop:cn},
$\tilde T$ is closed, so $\tilde t^*\in \tilde T$.  As we have just seen,
$|(u^*i)\tilde t^*|=1$, so $|\tilde t(ui)|=|(\tilde t(ui))^*|=1$.
\end{proof}

\section{Example}\label{sec:example}

We now provide an example to illustrate the definition above.
We let $X=\{1, 2, 3\}$, with basepoint $1$, and
let $T=\{1, t\}$ be the unique scheme on $X$ with $2$ elements.  That is, the adjacency matrix for $1$
is the identity, and the adjacency matrix for $t$ is
\[\sigma_t=\left(\begin{matrix} 0 & 1 & 1 \\ 1 & 0 & 1\\ 1 & 1 & 0 \end{matrix}\right).\]
Then $\tau=\tau_T=\{1, t\}$.
Let $Y$ be the set $\{a, b, c, d\}$, with basepoint $a$, and let $U$ be the thin scheme on $Y$
corresponding to the
group $\mathbb Z/4$.  Thus, the adjacency matrices for the elements of $Y$ are
\[\sigma_0=\left(\begin{matrix} 1 & 0 & 0 & 0 \\ 0 & 1 & 0 & 0\\ 0 & 0 & 1 & 0\\ 0 & 0 & 0 & 1
\end{matrix}\right),
\sigma_1=\left(\begin{matrix} 0 & 0 & 0 & 1 \\ 1 & 0 & 0 & 0 \\ 0 & 1 & 0 & 0\\ 0 & 0 & 1 & 0
\end{matrix}\right),
\sigma_2=\left(\begin{matrix} 0 & 0 & 1 & 0 \\ 0 & 0 & 0 & 1 \\ 1 & 0 & 0 & 0\\ 0 & 1 & 0 & 0
\end{matrix}\right),
\sigma_3=\left(\begin{matrix} 0 & 1 & 0 & 0 \\ 0 & 0 & 1 & 0 \\ 0 & 0 & 0 & 1\\ 1 & 0 & 0 & 0
\end{matrix}\right).\]
For each $y\in Y$, we let $T_y=T$ and $\alpha^y$ be the identity; then $\zeta^y=(T_y,\alpha^y)$.

Now, we define $\zeta_a^b\in \HC(T_a,T_b)$ as follows.  The normal closed subsets $(T_a)_{\zeta_a^b}$
and $(T_b)_{\zeta_a^b}$ are both $T$, and $\tilde\zeta_a^b$ is the identity on $T\dm T$.  We
define the seven morphisms $\zeta_b^c, \zeta_c^d, \zeta_d^a, \zeta_b^a, \zeta_c^b, \zeta_d^c, \zeta_a^d$
in the same way.  Next, we define $\zeta_a^c\in \HC(T_a, T_b)$ as follows.  The normal closed subsets
$(T_a)_{\zeta_a^c}$ and $(T_c)_{\zeta_a^c}$ are both $\{1\}$, and $\tilde\zeta_a^c$ is the identity
on $T\dm \{1\}\iso T$.  We define the three morphisms $\zeta_b^d, \zeta_c^a, \zeta_d^b$ in
the same way.  Finally, we define $\zeta_y^y$ to be the identity for each $y\in \{a, b, c, d\}$.
The first four conditions of Definition \ref{defn:action} are immediate.  The fifth condition is also easy.
Indeed, if either $(y_1, y_2)$ or $(y_2, y_3)$ belong to scheme elements $1$ or $3$, then the composition
$\zeta_{y_1}^{y_2}\zeta_{y_2}^{y_3}$ is the greatest element in the partial ordering of
$\HC(T_{y_1}, T_{y_3})$.  If $(y_1, y_2)$ and $(y_2, y_3)$
both belong to $2$, then $y_1=y_3$, and $\zeta_{y_1}^{y_2}\zeta_{y_2}^{y_3}$ is the identity
in $\HC(T_{y_1}, T_{y_3})$, and so is equal to $\zeta_{y_1}^{y_3}$.  Finally, if $y_1=y_2$
or $y_2=y_3$, then $\zeta_{y_1}^{y_2}\zeta_{y_2}^{y_3}=\zeta_{y_1}^{y_3}$.

Let $S=U\ltimes_\zeta T$.  Then $X\times Y$ has $12$ elements, and $S$
has $6$ elements: \[[0,1],[0,t],[1,1],[2,1],[2,t],[3,1].\]
Indeed, $[1,t]=[1,1]$ and $[3,t]=[3,1]$.  The elements $[0,1]$ and $[2,1]$ are both thin.
The valencies of $[0,t]$ and $[2,t]$ are both $2$, and the valencies of $[1,1]$ and $[3,1]$
are both $3$.  The elements of valency $1$ and $2$ are symmetric, while $[1,1]=[3,1]^*$.
Using Hanaki and Miyamoto's classification \cite{HM}, we see that there is only one
scheme of order $12$ satisfying these conditions, scheme No. 34.  Indeed, we can give
an explicit isomorphism between scheme No. 34 and our semidirect product.  On underlying
sets, this isomorphism would take the elements in the first row of the matrix below to the elements
directly below them in $X\times Y$:
\[\left(\begin{array}{cccccccccccc} 1 & 2 & 3 & 4 & 5 & 6 & 7 & 8 & 9 & 10 & 11 & 12\\
(a,1) & (c,1) & (a,2) & (a,3) & (c,2) & (c,3) & (b,1) & (b,2) & (b,3) & (d,1) & (d,2) & (d,3)
\end{array}\right)\]
On scheme elements, the isomorphism can be described by the matrix below:
\[\left(\begin{array}{cccccc} 0 & 1 & 2 & 3 & 4 & 5\\
\left[0,1\right] & [2,1] & [0,t] & [2,t] & [1,1] & [3,1]
\end{array}\right)\]

We next show that $S$ cannot be obtained
as a generalized semidirect product in the sense of Bang and Song \cite{BS}, and for the remainder
of this section, we will assume the reader is familiar with that paper.  First, we observe
that the scheme $S$ has exactly three proper non-trivial closed subsets:
\begin{itemize}
	\item $K_1=\{[0,1],[2,1]\}$
	\item $K_2=\{[0,1],[0,t]\}$
	\item $K_3=\{[0,1],[0,t],[2,1],[2,t]\}$
\end{itemize}
To see this, note first that $[1,1]$ and $[3,1]$ both generate $S$,
so any proper closed subset can only contain the other four elements.  Moreover, of the three elements
$[0,t]$, $[2,1]$, and $[2,t]$, each is in the complex product of the other two.

If $i=1$ or $i=2$, then $K_i$ has only two elements, so that $\text{Aut}(K_i)$ is trivial, and
$K_i$ has no nontrivial proper closed subsets.  It follows that the only generalized semidirect
products of $K_i$ by another scheme are the ordinary product and the wreath product.
Suppose $H$ were a scheme whose elements had valencies $h_1=1, h_2, \ldots, h_m$.
Then the valencies
of the wreath product of $H$ and $K_1$ would be $1, 1, 2h_2, 2h_3, \ldots, 2h_m$, while
those of the wreath product of $H$ and $K_2$ would be $1, 2, 3h_2, 3h_3, \ldots, 3h_m$.
In either case, it is impossible to obtain valencies $1, 1, 2, 2, 3, 3$.  On the other hand,
the valencies of the direct product $H\times K_2$ would be $h_1, h_2, \ldots, h_m, 2h_1, 2h_2, \ldots, 2h_m$.
Again, we cannot obtain $1, 1, 2, 2, 3, 3$.  We could obtain this sequence
of valencies from $H\times K_1$ if the valencies of $H$ were $1, 2, 3$.  However,
such a scheme must clearly be symmetric, so $H\times K_1$ would be symmetric.
But $S$ is not symmetric, since the two elements of valency $3$ are conjugates.

Now, we claim $\text{Aut}(K_3)$ is also trivial.  Indeed, any automorphism must send $[0,1]$ to itself,
since this represents the identity.  Since $[2,1]$ is the only other thin element, any automorphism
must send $[2,1]$ to itself.  Thus, the only possible nontrivial automorphism would transpose
$[0,t]$ and $[2,t]$.  However, if we let $p=[0,t]$ and $q=[2,t]$, then $a_{ppp}=1$, while
$a_{qqq}=0$.   Thus, transposing $[0,t]$ and $[2,t]$ does not determine an automorphism of $K_3$.

Thus, any semidirect product in the sense of Bang, Hirasaka, and Song \cite{BHS} of $K_3$
by another scheme must be a direct product.
In order that such a semidirect product have order $12$, the other scheme would have
to be the unique scheme of order $2$, since the valency of $K_3$ is $6$.
Since $K_3$ and the scheme of order $2$ are both symmetric, such a direct product would
be symmetric.  Since
the generalized semidirect product of Bang and Song \cite{BS} is a fusion of the semidirect product
in the sense of Bang, Hirasaka and Song \cite{BHS} by Theorem 2.1 in the former paper,
and any fusion of a symmetric scheme is symmetric, we cannot obtain $S$ by taking a
generalized semidirect product of $K_3$ by a scheme of order $2$.

\section{Actions obtained from semidirect products}\label{sec:action}

In this section, we suppose given two based schemes $T$ on $X$ and $U$ on $Y$.  We suppose
$S$ is a based scheme on a set $Z$ (with basepoint $z_*$), equipped with a closed subset $\tilde T\subseteq S$ such that
there is a based isomorphism from the subscheme of $S$ defined by $z_*\tilde T$ to the scheme $T$.
We also assume that there is a based morphism $i:U\to S$ such that the composition
$i\pi:U\to S\dm \tilde T$ is an isomorphism, where $\pi:S\to S\dm \tilde T$ denotes the natural morphism.
Finally, we assume Condition \ref{con:hyp}:  for any $u\in U$ and $t\in \tilde T$, we have $|t(ui)|=1$.
Our goal in this final section is to show that one can construct an action $\zeta$ of $U$ on $T$
such that $U\ltimes_\zeta T$ is isomorphic to $S$.

\begin{lem}\label{lem:TTnormal} Given the conditions above, $\tilde T$ is normal.
\end{lem}

\begin{proof}
Let $z=yi$ for some $y\in Y$.
To prove the claim, it suffices to show that for any $s\in S$, we have $\tilde Ts\subseteq s\tilde T$.  
We first assume that $s=ui$ for some $u\in U$.  If $z'\in z\tilde T(ui)$, then $z'\pi\in z\pi(ui\pi)=yi\pi(ui\pi)$.
Since $i\pi$ is an isomorphism, $yi\pi(ui\pi)=(yu)i\pi\subseteq ((yi)(ui))\pi$.  Thus, $z'\pi\in ((yi)(ui))\pi=z(ui)\pi$.
This implies $z'\in z(ui)\tilde T$.  Since $z'$ was arbitarily chosen from $z\tilde T(ui)$, we have
$z\tilde T(ui)\subseteq z(ui)\tilde T$, whence $\tilde T(ui)\subseteq (ui)\tilde T$.

Now, given any $s\in S$, we have $s\pi=u(i\pi)$ for some $u\in U$, since $i\pi$ is an isomorphism.
That is, $s^{\tilde T}=(ui)^{\tilde T}$, so $s\in \tilde T (ui) \tilde T$.  Since $\tilde T(ui)\subseteq (ui)\tilde T$,
we find $s\in (ui)\tilde T$.  Thus, for some $t\in \tilde T$, $s\in (ui)t$.  By Condition \ref{con:hyp}, $|(ui)t|=1$,
so $(ui)t=\{s\}$.  This implies $(ui)\tilde T=s\tilde T$.  Now, since $\tilde T(ui)\subseteq (ui)\tilde T$, we obtain
$\tilde T s\subseteq \tilde T (ui) \tilde T\subseteq (ui)\tilde T=s\tilde T.$
\end{proof}

Next, we wish to define an action $\zeta$ of $U$ on $T$.  For this, we will need a $\tau$-scheme
$\zeta_y=(T_y,\al^y)$ on $X$ for each $y\in Y$, where $\tau=\tau_T$.
For each $y\in Y$, we first let $\tilde T_y$ denote the
subscheme of $S$ defined by $(yi)\tilde T$, with basepoint $yi$.  That is, $\tilde T_y$ consists of elements
$t\cap (yi\tilde T\times yi\tilde T)$, for $t\in \tilde T$.  Let $\de^y:\tilde T_y\to \tilde T$ be
the bijection taking $t\cap (yi\tilde T\times yi\tilde T)$ to $t$.

We had assumed at the beginning of this section that the subscheme of $S$ defined by
$z_*\tilde T$ is isomorphic
to $T$.  Thus, since $z_*=y_*i$, we may choose a based isomorphism $\ga:T\to \tilde T_{y_*}$. 
(Thus, $\ga_T$ is a bijection from $T$ to $\tilde T_{y_*}$, and $\ga_X$ is a bijection from
$X$ to $(y_*i)\tilde T$ taking $x_*$ to $y_*i$.)  For each $y\in Y$, we have $|(yi)\tilde T|=|z_*\tilde T|=|X|$,
so we may 
choose a set of based bijections $\ga^y_X:X\to (yi)\tilde T$, one for each $y\in Y$;
we choose $\ga^{y_*}_X$ to be $\ga_X$.

\begin{defn}\label{def:TY}
Let $T_y$ be the unique scheme on $X$ defined by the requirement that $\ga^y_X:X\to (yi)\tilde T$
determines a based isomorphism of schemes $\ga^y:T_y\to \tilde T_y$.
Let $\al^y:\tau_T\to T_y$ be the bijection defined by $(t\al^y)\ga^y_{T_y}\de^y=t\ga_T\de^{y_*}$
for each $t\in \tau_T=T$.  Let $\zeta_y=(T_y,\al^y)$.
\end{defn}

Thus, the elements of $T_y$ are the preimages of elements in $\tilde T_y$ under the product
bijection \[\ga^y_X\times \ga^y_X:X\times X\to (yi)\tilde T\times (yi)\tilde T.\]
Since $\ga:T\to \tilde T_{y_*}$ is already an isomorphism of schemes
and $\ga^{y_*}_X=\ga_X$, it then follows that $T_{y_*}=T$.
Also, since $\ga_T=\ga^{y_*}_{T_{y_*}}$, it follows that $\al^{y_*}$ is the identity.
The following commutative diagram, in which $y$ is an arbitrary element in $Y$ and $y_*$
is the basepoint, may be helpful in keeping these definitions straight:

\begin{equation}\label{diagram}
\xymatrix{&& {T_{y}}\ar[rr]^-{\ga^{y}_{T_{y}}} && {\tilde T_{y}}\ar[drr]^-{\de^{y}} & & \\
{\tau}\ar[urr]^-{\al^{y}}\ar@{=}^{\al^{y_*}}[rr] &
& {T_{y_*}}\ar[rr]^-{\ga_T=\ga_{T_{y_*}}^{y_*}} & &
{\tilde T_{y_*}}\ar[rr]^{\de^{y_*}} & & {\tilde T} \\ }
\end{equation}

We next require morphisms $\zeta_{y_1}^{y_2}\in \HC(T_{y_1},T_{y_2})$ for each pair $y_1, y_2\in Y$.
Recall that such a morphism consists of normal closed subsets $T'_{y_1y_2}$ in $T_{y_1}$ and $T''_{y_1y_2}$
in $T_{y_2}$, together with a based isomorphism of schemes $\tilde\zeta_{y_1}^{y_2}$ from
$T_{y_1}\dm T'_{y_1y_2}$ to $T_{y_2}\dm T''_{y_1y_2}$.  The following subsets of $\tilde T$
will be useful in defining $T'_{y_1y_2}$ and $T''_{y_1y_2}$.

\begin{defn}\label{def:T'T''} For a given $u\in U$, we define
\[\tilde T'_u=\{t\in \tilde T:t(ui)=\{ui\}\}\text{ and }
\tilde T''_u=\{t\in \tilde T:(ui)t=\{ui\}\}.\]
Let $\tau_u'=\{t\in \tau:t\ga_T\de^{y_*}\in \tilde T'_u\}$,
and let $\tau_u''=\{t\in \tau:t\ga_T\de^{y_*}\in \tilde T''_u\}$.
\end{defn}

\begin{lem}\label{lem:closed} The sets $\tilde T'_u$ and $\tilde T''_u$ are closed in $\tilde T$.
\end{lem}

\begin{proof}
Suppose $p, q\in \tilde T'_u$.  Then $p(ui)=\{ui\}$ and $q(ui)=\{ui\}$.  Since
$ui\in p^*p(ui)=p^*(ui)$, and $p^*(ui)$ consists of precisely one element (by Condition \ref{con:hyp}),
it follows that $p^*(ui)=\{ui\}$.  Therefore, $p^*q(ui)=\{ui\}$, so $p^*q\subseteq \tilde T'_u$.
Thus, $\tilde T'_u$ is closed, and by a similar argument, $\tilde T''_u$ is closed as well.
\end{proof}

\begin{defn}\label{def:T'T''2}  Suppose
$(y_1,y_2)\in u$.  Let $T'_{y_1y_2}=\tau'_u\al^{y_1}$ and $T''_{y_1y_2}=\tau''_u\al^{y_2}$.
Also let $\tilde T'_{y_1y_2}=T'_{y_1y_2}\ga_{T_{y_1}}^{y_1}$, and let $\tilde T''_{y_1y_2}=
T''_{y_1y_2}\ga_{T_{y_2}}^{y_2}$.  
\end{defn}

\begin{rem}\label{rem:delta}
By the commutativity of Diagram (\ref{diagram}), and the fact that $\al^y$ and $\ga^{y}_{T_y}$
are bijections for each $y\in Y$, it follows that if $(y_1,y_2)\in u$, then 
$\de^{y_1}$ defines a bijection from $\tilde T'_{y_1y_2}$ to $\tilde T_u'$
and $\de^{y_2}$ defines a bijection from $\tilde T''_{y_1y_2}$ to $\tilde T_u''$.
Also, from the way they are defined,
it is clear that $\de^{y_1}$ and $\de^{y_2}$ preserve the involution and the
complex product on subsets of $\tilde T'_{y_1y_2}$ and $\tilde T''_{y_1y_2}$.
Since $\ga^{y_1}$ and $\ga^{y_2}$ are isomorphisms of schemes, it follows now from
Lemma \ref{lem:closed} that
the sets $T'_{y_1y_2}$ and $T''_{y_1y_2}$ are closed subsets of $T_{y_1}$ and $T_{y_2}$ respectively.
\end{rem}

In order to show that $T'_{y_1y_2}$ and $T''_{y_1y_2}$ are normal in $T_{y_1}$ and $T_{y_2}$,
we will need to show that $\tilde T'_u$ and $\tilde T''_u$ are normal in $\tilde T$, which
is more difficult than proving that they are closed.  We
will require a few preliminary lemmas.  First, we recall that if $u$ is an element
of a scheme $U$ on a set $Y$, then the valency of $u$, denoted $n_u$, is equal to $a_{uu^*1}$.  For any
$y\in Y$, $n_u$ is the number of elements in $yu$.  Similarly, if $R\subseteq U$
is any subset, then $n_R=\sum_{u\in R} n_u$ is the number of elements in $yR$.

\begin{lem}\label{lem:valence}  Suppose $u\in U$.  Then $n_{ui}=n_u\cdot n_{\tilde T''_u}.$
\end{lem}

\begin{proof} We will show that 
\[(y_1i)(ui)=\bigsqcup_{y\in y_1u} (yi)\tilde T''_{u}.\]
That is, the set $(y_1i)(ui)$ decomposes as a disjoint union of $n_u$ sets,
each having $n_{\tilde T''_{u}}$ elements.  (Note that the union is disjoint since $i\pi$ is an injective
and $\tilde T''_u\subseteq \tilde T$.)
If $y\in y_1u$, then $yi\in (y_1i)(ui)$.  But $(ui)t=\{ui\}$ for any $t\in \tilde T''_{u}$,
so $(yi)\tilde T''_{u}\subseteq (y_1i)(ui)$.  Thus, \[\bigsqcup_{y\in y_1u} (yi)\tilde T''_{u}\subseteq (y_1i)(ui).\]

Now, suppose $z\in (y_1i)(ui)$.  Since $i\pi$ is an isomorphism, we have $z\pi=yi\pi$ for some $y\in Y$, and
also $z\pi\in (y_1i\pi)(ui\pi)=(y_1u)i\pi$.  Thus, $y\in y_1u$, so $yi\in (y_1i)(ui)$.  Now,
since $z\pi=yi\pi$, $z\in (yi)t\subseteq (y_1i)(ui)t$ for some $t\in \tilde T$.
But we supposed to begin with that $z\in (y_1i)(ui)$.  Thus,
we must have $ui\in (ui)t$, whence $(ui)t=\{ui\}$ by Condition \ref{con:hyp}.  Therefore,
$t\in \tilde T''_{u}$.  Thus, $z\in (yi)\tilde T''_{u}$ for some $y\in y_1u$.   Since $z$
was an arbitrary element in $(y_1i)(ui)$, we have
\[(y_1i)(ui)\subseteq \bigsqcup_{y\in y_1u} (yi)\tilde T''_{u}.\]
\end{proof}

\begin{lem}\label{lem:intersect} Suppose $(y_1,y_2)\in u$.  If
$z_1\in (y_1i)\tilde T$, then
$z_1(ui)\cap (y_2i)\tilde T=z_2\tilde T''_{u}$
for some $z_2\in (y_2i)\tilde T$.  Moreover, if $z_1'\in z_1 \tilde T'_{u}$, then
$z_1'(ui)\cap (y_2i)\tilde T=z_1(ui)\cap (y_2i)\tilde T$.
\end{lem}

\begin{proof} 
First, since $z_1\tilde T=y_1i\tilde T$, the cosets of $\tilde T$ which contain elements in $z_1(ui)$ are
the same as the cosets of $\tilde T$ containing elements in $(y_1u)i$, and there are $n_u$ such cosets.
By Lemma \ref{lem:valence}, $z_1(ui)$ contains $n_u\cdot n_{\tilde T''_{u}}$ elements.
Moreover, $z_1(ui)$ decomposes into cosets of $\tilde T''_{u}$, each of which contains
$n_{\tilde T''_{u}}$ elements.  Thus, there must be exactly $n_u$ cosets of $\tilde T''_{u}$
in $z_1(ui)$, and each of these cosets is contained in a different coset of $\tilde T$
since $\tilde T''_{u}\subseteq \tilde T$.  Thus for any $y_2\in y_1u$, 
$z_1(ui)\cap (y_2i)\tilde T$ must coincide with
one of the cosets of $\tilde T''_{u}$.  If we choose an element $z_2$ in this coset,
then $z_1(ui)\cap (y_2i)\tilde T=z_2\tilde T''_{u}$.
For the last statement, if $z_1'\in z_1\tilde t$ for some $\tilde t\in \tilde T'_u$, then since
$\tilde t(ui)=\{ui\}$, we have $z_1(ui)=z_1'(ui)$.
\end{proof}

\begin{lem}\label{lem:s}
Suppose given $u\in U, t\in \tilde T$ and $s\in (ui)t$.  If $\tilde t\in \tilde T''_u$, then $s\tilde t=\{s\}$.
\end{lem}

\begin{proof}
Choose $(y_1, y_2)\in u$, and let $z_1=y_1i$ and $z_2=y_2i$, so $(z_1,z_2)\in ui$.  Since $s\in (ui)t$,
we have $ui\in st^*$, so we can find $z_3\in z_1s\cap z_2t$.  Now, if $r\in s\tilde t$, then
$s\in r\tilde t^*$, so we can find $z_4\in z_1r\cap z_3\tilde t$.  By Lemma \ref{lem:intersect} (applied to
$u^*$), $z_3(ui)^*\cap z_1\tilde T$ is nonempty, so we may choose $z_5\in z_3(ui)^*\cap z_1t'$
for some $t'\in \tilde T$.  Since $z_3\in z_1s$ and $z_5\in z_3(ui)^*\cap z_1t'$, we have
$s\in t'(ui)$, so by Condition \ref{con:hyp}, $t'(ui)=\{s\}$.
Since $\tilde t\in \tilde T''_u$, we have $(ui)\tilde t=\{ui\}$, so we must have
$z_4\in z_5(ui)$.  But then $z_4\in z_1t'(ui)=z_1s$ and $z_4\in z_1r$,
so $r=s$.  Since $s\tilde t$ is nonempty and can only contain $s$, it must be equal to $\{s\}$.
\end{proof}

\begin{prop}\label{prop:normal} For each $u\in U$, $\tilde T''_u$ is normal in $\tilde T$.
\end{prop}

\begin{proof} It will suffice to show
that $\tilde T''_u \tilde t\subseteq \tilde t\tilde T''_u$ for any $\tilde t\in \tilde T$.
Note that if $\tilde t\in \tilde T$, then
$\tilde t^*\in \tilde T$ since $\tilde T$ is closed.  If $r\in \tilde T''_u\tilde t$, then $r\in \tilde t''\tilde t$ for
some $\tilde t''\in \tilde T''_u$.
Then we can choose $z_1, z_2, z_3\in Z$ such that $z_3\in z_1r\cap z_2\tilde t$ and $z_2\in z_1\tilde t''$,
and we may assume $z_3=y_3i$ for some $y_3\in Y$.  Now, choose $y\in y_3u^*$, so $y_3\in yu$.
Let $s\in S$ be the element containing $(yi,z_2)$.  Then $s\in (ui)\tilde t^*$.  By Lemma \ref{lem:s},
$s(\tilde t'')^*=\{s\}$, so $z_1\in (yi)s$.  Since $s\in (ui)\tilde t^*$, we may find $z_4\in (yi)(ui)\cap z_1\tilde t$.
Note that $z_1$, $z_2$, $z_3$ and $z_4$ are all in the same coset of $\tilde T$, so $(z_3, z_4)\in \tilde t'$
for some $\tilde t'\in \tilde T$.  But then $ui\in (ui)(\tilde t')^*$, so $(ui)(\tilde t')^*=\{ui\}$,
and $(\tilde t')^*\in \tilde T''_u$.
Finally, $r\in \tilde t(\tilde t')^*\subseteq \tilde t \tilde T''_u$,
as we see by considering $z_1$, $z_3$ and $z_4$.  So, $\tilde T''_u \tilde t
\subseteq \tilde t T''_u$, as needed.
\end{proof}

\begin{cor} $\tilde T'_u$ is normal in $\tilde T$.
\end{cor}

\begin{proof}
It is easy to check that $(\tilde T'_{u})^*=\tilde T''_{u^*}$, and by Lemma \ref{lem:closed},
$(\tilde T'_u)^*=\tilde T'_u$, so the normality of $\tilde T'_u$ follows from Proposition \ref{prop:normal}
applied to $u^*$.
\end{proof}

The following corollary now follows by Remark \ref{rem:delta}, and since $\ga^{y_1}$ and
$\ga^{y_2}$ are isomorphisms of schemes.

\begin{cor}\label{cor:normal}
For any $y_1, y_2\in Y$, the sets $T'_{y_1y_2}$ and $T''_{y_1y_2}$ are normal in $T_{y_1}$
and $T_{y_2}$.
\end{cor}

Next, we need to define for each pair $y_1, y_2\in Y$ a based isomorphism of
schemes $\tilde\zeta_{y_1}^{y_2}$ from $T_{y_1}\dm T'_{y_1y_2}$ to $T_{y_2}\dm T''_{y_1y_2}$.
Using the isomorphisms $\ga^{y_1}$ and $\ga^{y_2}$, it will suffice to define a based
isomorphism $\tilde\xi_{y_1}^{y_2}$ from
$\tilde T_{y_1}\dm \tilde T'_{y_1y_2}$ to $\tilde T_{y_2}\dm \tilde T''_{y_1y_2}$.
We suppose $u\in U$ is the element containing $(y_1, y_2)$.

\begin{defn}\label{def:xi}
Let \[\tilde \xi_{y_1}^{y_2}:(y_1i)\tilde T/\tilde T'_{y_1y_2}\to (y_2i)\tilde T/\tilde T''_{y_1y_2}\]
be the function which takes the coset $z_1\tilde T'_{y_1y_2}$ to the coset $z_2\tilde T''_{y_1y_2}$,
where $z_2$ is any element in $z_1(ui)\cap (y_2i)\tilde T$.
\end{defn}

\begin{lem} The function $\tilde \xi_{y_1}^{y_2}$ is well-defined.
\end{lem}

\begin{proof}
If $z_1\in (y_1i)\tilde T$, then by Lemma \ref{lem:intersect}, there is a $z_2\in (y_2i)\tilde T$ such that
$z_1(ui)\cap (y_2i)\tilde T=z_2\tilde T''_u$, and the coset $z_2\tilde T''_u$ only
depends on the coset of $\tilde T'_u$ containing $z_1$.  The lemma will follow if we can show
that the coset of $\tilde T'_u$
containing $z_1$ is the same as the coset of $\tilde T'_{y_1y_2}$ containing $z_1$, and that
the coset of $\tilde T''_u$ containing $z_2$ is the same as the coset of $\tilde T''_{y_1y_2}$ containing
$z_2$.  Indeed, restricting elements of $\tilde T'_u$
and $\tilde T''_u$ to $(y_1i)\tilde T$ and $(y_2i)\tilde T$ corresponds to taking the preimage under
$\delta^{y_1}$ and $\delta^{y_2}$.  Thus, the restrictions of elements of $\tilde T'_u$ and
$\tilde T''_u$ to $(y_1i)\tilde T$ and $(y_2i)\tilde T$ coincide with $\tilde T'_{y_1y_2}$
and $\tilde T''_{y_1y_2}$ by Definitions \ref{def:T'T''} and \ref{def:T'T''2}.

\end{proof}

\begin{lem}\label{lem:xi}
The function $\tilde\xi_{y_1}^{y_2}$ determines a based morphism of schemes from
$\tilde T_{y_1}\dm \tilde T'_{y_1y_2}$ to $\tilde T_{y_2}\dm \tilde T''_{y_1y_2}$.  Moreover, 
for $v\in \tilde T'_{y_1}$, we have
$v^{\tilde T'_{y_1y_2}}\tilde\xi_{y_1}^{y_2}=(v')^{\tilde T''_{y_1y_2}}$ if and only if
$(v\de^{y_1})(ui)=(ui)(v'\de^{y_2})$.
\end{lem}

\begin{proof}
Since $y_2i\in (y_1i)(ui)$ it follows that $\tilde\xi_{y_1}^{y_2}$ takes the coset of $\tilde T'_{y_1y_2}$
containing $y_1i$ to the coset of $\tilde T''_{y_1y_2}$ containing $y_2i$, so $\tilde\xi_{y_1}^{y_2}$ is based.
Suppose that $z_1, z_1'\in (y_1i)\tilde T$, and let $v\in \tilde T_{y_1}$ be the scheme element
containing $(z_1, z_1')$.  Choose $z_2, z_2'$ as above such that
$z_1(ui)\cap (y_2i)\tilde T=z_2\tilde T''_u$ and $z_1'(ui)\cap (y_2i)\tilde T=z_2'\tilde T''_u$,
and let $v'\in \tilde T_{y_2}$ be the scheme element containing $(z_2, z_2')$.  Then since $z_2\in z_1(ui)$
and $z_2'\in z_1'(ui)$, we have $(v\de^{y_1})(ui)\cap (ui)(v'\de^{y_2})\neq \emptyset$.
By Condition \ref{con:hyp}, $(v\de^{y_1})(ui)=(ui)(v'\de^{y_2})$.  Note that if $(ui)r=(ui)r'$
for two elements $r, r'\in \tilde T$, then $ui\in (ui)r'r^*$, so for some $p\in r'r^*\subseteq \tilde T$,
$ui\in (ui)p$, whence $(ui)p=\{ui\}$, so $p\in \tilde T''_u$.  Since $r'\in pr$, we have
$(r)^{\tilde T''_u}=(r')^{\tilde T''_u}$.  Thus, the equation
$(v\de^{y_1})(ui)=(ui)(v'\de^{y_2})$ implies that $(v'\de^{y_2})^{\tilde T''_u}$
is uniquely determined by $v\de^{y_1}$.  By Remark \ref{rem:delta}, this
in turn implies that $(v')^{\tilde T''_{y_1y_2}}$ is uniquely determined
by $v$.  Similarly, if $r^{\tilde T'_u}=r'^{\tilde T'_u}$, then $r\in r't$ for some $t\in \tilde T'_u$
(since $\tilde T'_u$ is normal), so
$r(ui)=r'(ui)$.  Thus, $(v')^{\tilde T''_{y_1y_2}}$ is uniquely determined by
$v^{\tilde T'_{y_1y_2}}$, so $\tilde\xi_{y_1}^{y_2}$ determines a morphism of schemes.

For the second statement, if $v^{\tilde T'_{y_1y_2}}\tilde \xi_{y_1}^{y_2}
=(v')^{\tilde T''_{y_1y_2}}$, then we can find $(z_1, z_1')\in v$ and choose $(z_2, z_2')$
as above so that $(z_2, z_2')\in v'$.  But then, as shown above, we have
$(v\de^{y_1})(ui)=(ui)(v'\de^{y_2})$.  Conversely, if $(v\de^{y_1})(ui)=(ui)(v'\de^{y_2})$,
then again, we may choose $(z_1, z_1')\in v$.  We choose $z_2'$ in 
$z_1'(ui)\cap (y_2i)\tilde T$.  Then, $z_2'\in z_1(v\de^{y_1})(ui)=z_1(ui)(v'\de^{y_2})$,
so there must exist a $z_2\in z_1(ui)$ with $(z_2, z_2')\in (v'\de^{y_2})$.  Then
$z_2\in z_1(ui)\cap (y_2i)\tilde T$ since $(v'\de^{y_2})\in \tilde T$ and $z_2'\in (y_2i)\tilde T$.
But then, we see as above that $(v)^{\tilde T'_{y_1y_2}}\tilde\xi_{y_1}^{y_2}=(v')^{\tilde T''_{y_1y_2}}$.
\end{proof}

\begin{lem}\label{lem:xi2}
The morphism $\tilde\xi_{y_1}^{y_2}$ is an isomorphism, with inverse $\tilde\xi_{y_2}^{y_1}$.
\end{lem}

\begin{proof}
First, we need to see that $\tilde T'_{y_2y_1}=\tilde T''_{y_1y_2}$ and $\tilde T''_{y_2y_1}=\tilde T'_{y_1y_2}$.
Following Definition \ref{def:T'T''}, we have
\[\tilde T'_{u^*}=\{t\in \tilde T:t(u^*i)=\{u^*i\}\}
=\{t\in \tilde T:(ui)t^*=\{ui\}\}=(\tilde T''_u)^*.\]  Similarly, $\tilde T''_{u^*}=(\tilde T'_{u})^*$.
Since $\ga_T$ and $\de^{y^*}$ preserve the involution, it follows that $\tau'_{u^*}=(\tau''_u)^*$
and $\tau''_{u^*}=(\tau'_u)^*$, and then since $\al^{y_1}$ and $\ga^{y_1}_{T_{y_1}}$
preserve the involution
and $\tilde T'_{y_1y_2}$ is closed, we have
\[\tilde T''_{y_2y_1}=\tau''_{u^*}\al^{y_1}\ga^{y_1}_{T_{y_1}}
=(\tau'_u\al^{y_1}\ga^{y_1}_{T_{y_1}})^*=(\tilde T'_{y_1y_2})^*=\tilde T'_{y_1y_2}.\]
Similarly $\tilde T'_{y_2y_1}=\tilde T''_{y_1y_2}$.

Now, given $z_1\in (y_1i)\tilde T$, we have
$(z_1\tilde T'_{y_1y_2})\tilde \xi_{y_1}^{y_2}=z_2\tilde T''_{y_1y_2}$
if and only if $z_2\in z_1(ui)\cap (y_2i)\tilde T$.  But if $z_2\in z_1(ui)$, then $z_1\in z_2(u^*i)$,
so $z_1\in z_2(u^*i)\cap (y_1i)\tilde T$.  Thus, if
$(z_1\tilde T'_{y_1y_2})\tilde \xi_{y_1}^{y_2}=z_2\tilde T''_{y_1y_2}$,
then
$(z_2\tilde T'_{y_2y_1})\tilde \xi_{y_2}^{y_1}=z_1\tilde T''_{y_2y_1}.$
The converse follows by reversing the rolls of $y_1$ and $y_2$.  Therefore, $\tilde \xi_{y_2}^{y_1}$
is the inverse of $\tilde \xi_{y_1}^{y_2}$.
\end{proof}

As mentioned in the paragraph before Definition \ref{def:xi}, the isomorphisms $\tilde\xi_{y_1}^{y_2}$
will induce the needed isomorphisms $\tilde\zeta_{y_1}^{y_2}$, using the isomorphisms $\gamma^{y_1}$
and $\gamma^{y_2}$.  This is made precise in the following definition.

\begin{defn}\label{def:tildezeta}
Given $y_1, y_2\in Y$, let $\tilde\zeta_{y_1}^{y_2}$ be the based
isomorphism of schemes making the following diagram commute:

\[\xymatrix{
{T_{y_1}\dm T'_{y_1y_2}}\ar[r]^-{\overline{\ga^{y_1}}}\ar[d]^-{\tilde\zeta_{y_1}^{y_2}} &
{\tilde T_{y_1}\dm \tilde T'_{y_1y_2}}\ar[d]^-{\tilde\xi_{y_1}^{y_2}} \\
{T_{y_2}\dm T''_{y_1y_2}}\ar[r]^-{\overline{\ga_{y_2}}} &
{\tilde T_{y_2}\dm \tilde T''_{y_1y_2}}}\]
\end{defn}

The horizontal arrows in the above diagram are the isomorphisms of quotient
schemes induced by the isomorphisms $\ga^{y_1}$ and $\ga^{y_2}$.
(Recall that $\tilde T'_{y_1y_2}$ is defined simply as the image of $T'_{y_1y_2}$
under $\gamma^{y_1}$ and $\tilde T''_{y_1y_2}$ is defined as the image of $T''_{y_1y_2}$
under $\gamma^{y_2}$.)

Having defined $\tau$-schemes $\zeta_y$ for each $y\in Y$ and morphisms
$\zeta_{y_1}^{y_2}\in \HC(T_{y_1},T_{y_2})$ for each pair of elements $y_1, y_2\in Y$,
we must now prove that $\zeta$ determines an action of $U$ on $T$.  That is, we need
to verify the conditions of Definition \ref{defn:action}.  The following lemma will be needed.

\begin{lem}\label{lem:valency2} Suppose $p, q, r\in U$, with $r\in pq$, and suppose $z_1, z_3\in Z$,
with $z_3\in z_1(ri)$.  Let $y_1$ and $y_3$ be the elements in $Y$ such that $z_1\tilde T=(y_1i)\tilde T$
and $z_3\tilde T=(y_3i)\tilde T$.  If $y_2\in y_1p\cap y_3q^*$, there exists an element
$z_2\in z_1(pi)\cap z_3(qi)^*$ such that $z_2\tilde T=(y_2i)\tilde T$.
\end{lem}

\begin{proof}
First, we show that for each $y_2\in y_1p\cap y_3q^*$, the number of elements
$z_2\in z_1(pi)\cap z_3(qi)^*$ such that $z_2\tilde T=(y_2i)\tilde T$
is either $0$ or $n_{\tilde T''_p\cap \tilde T'_q}$.  To see this,
it suffices to show that if $z_1(pi)\cap z_3(qi)^*\cap (y_2i)\tilde T$ contains an element $z_2$,
then it coincides with the coset of $\tilde T''_p\cap \tilde T'_q$ containing $z_2$.  Suppose
$t\in \tilde T''_p\cap \tilde T'_q$.  Then $t^*\in \tilde T'_q$, so
we have $z_2t\subseteq z_1(pi)t=z_1(pi)$ and
$z_2t\subseteq z_3(qi)^*t=z_3(t^*(qi))^*=z_3(qi)^*$.  Since
$\tilde T''_p\cap \tilde T'_q\subseteq \tilde T$ and $z_2\in (y_2i)\tilde T$,
we have $z_2t\subseteq z_2\tilde T=(y_2i)\tilde T$.  Thus, every element in the coset of
$\tilde T''_p\cap \tilde T'_q$ containing $z_2$ is in $z_1(pi)\cap z_3(qi)^*\cap (y_2i)\tilde T$.
Conversely, if $z\in z_1(pi)\cap z_3(qi)^*\cap (y_2i)\tilde T$,
then let $t\in \tilde T$ be the element containing $(z_2,z)$.  
Considering the elements $z_1$, $z_2$ and $z$, we see that $pi\in (pi)t$, so 
$(pi)t=\{pi\}$, whence $t\in \tilde T''_p$;
considering $z_2$, $z$, and $z_3$, we see that $qi\in t(qi)$, so $t(qi)=\{qi\}$, whence
$t\in \tilde T'_q$.  Thus, $z$ is in the coset of $\tilde T''_p\cap \tilde T'_q$ containing $z_2$.

The argument in the previous paragraph
applies when we choose $z_1=y_1i$ and $z_3=y_3i$, but in this case,
$z_1(pi)\cap z_3(qi)^*\cap (y_2i)\tilde T$ always contains $y_2i$.  Thus,
for each $y_2\in y_1p\cap y_3q^*$,
there are $n_{\tilde T''_p\cap \tilde T'_q}$ elements $z_2\in (y_1i)(pi)\cap (y_3i)(qi)^*$ with
$z_2\tilde T=(y_2i)\tilde T$.  Of course, if $y_2\notin y_1p\cap y_3q^*$, then
$(y_1i)(pi)\cap (y_3i)(qi)^*\cap (y_2i)\tilde T$ is empty since $i\pi$ is an isomorphism.
Since the cardinality of $y_1p\cap y_3q^*$ is $a_{pqr}$, 
and the cardinality of $z_1(pi)\cap z_3(qi)^*$ is $a_{(pi)(qi)(ri)}$, we have
shown that
\begin{equation}\label{eqn:apiqiri}
a_{(pi)(qi)(ri)}=a_{pqr}n_{\tilde T''_p\cap \tilde T'_q}.
\end{equation}

Returning to the general case, where $z_1\in (y_1i)\tilde T$, $z_3\in (y_3i)\tilde T$,
and $z_3\in z_1(ri)$,
the cardinality $a_{(pi)(qi)(ri)}$ of $z_1(pi)\cap z_3(qi)^*$ is equal to the product of
$n_{\tilde T''_p\cap \tilde T'_q}$ with the number of elements $y_2$ in $y_1p\cap y_3q^*$ for which
$z_1(pi)\cap z_3(qi)^*\cap (y_2i)\tilde T$ is nonempty.  By Equation 
(\ref{eqn:apiqiri}), there must be $a_{pqr}$ such elements $y_2$.  Since $|y_1p\cap y_3q^*|=a_{pqr}$,
it must be the case that for every $y_2\in y_1p\cap y_3q^*$, we have $z_1(pi)\cap z_3(qi)^*\cap (y_2i)\tilde T$
is nonempty.
\end{proof}

\begin{thm}\label{prop:actin} Definitions \ref{def:TY}, \ref{def:T'T''2}, and \ref{def:tildezeta}
determine an action $\zeta$ of $U$ on $T$.
\end{thm}

\begin{proof} We must verify the five conditions of Definition \ref{defn:action}.

\medskip

\noindent {\bf Condition (1)}

We have seen in the paragraph following Definition \ref{def:TY} that condition (1) is satisfied.

\medskip

\noindent {\bf Condition (2)}

Suppose $y\in Y$, and let $1_Y\in U$
denote the scheme element containing $(y,y)$.  Then since $1_Yi=1_Z\in S$, it follows immediately
from Definition \ref{def:T'T''} that $\tilde T'_{1_Y}=\tilde T''_{1_Y}=\{1_Z\}$, and, since $\gamma_T$
is an isomorphism, and $\delta^{y}$ is a bijection taking $1_{(yi)\tilde T}$ to $1_Z$,
$\tau'_{1_Y}=\tau''_{1_Y}=\{1\}$.  Now it follows from Definition \ref{def:T'T''2} that
$T'_{yy}=T''_{yy}=\{1_X\}$.  Now, if $z_1\in (yi)\tilde T$, then clearly, $z_1(1_Z)\cap (yi)\tilde T=z_1(1_Z)$,
so $\tilde\xi_{y}^{y}$ is the identity.  By Definition \ref{def:tildezeta}, $\tilde\zeta_{y}^{y}$ is then
also the identity.  Thus, $\zeta_y^y$ is the identity morphism in $\HC(T_y,T_y)$, as needed for Condition (2).

\medskip

\noindent {\bf Condition (3)}

If $y_1, y_2\in Y$, then by Lemma \ref{lem:xi2}, $\tilde\xi_{y_2}^{y_1}=(\tilde\xi_{y_1}^{y_2})^{-1}$, 
so by Definition \ref{def:tildezeta}, $\tilde\zeta_{y_2}^{y_1}=(\tilde \zeta_{y_1}^{y_2})^{-1}$.
Therefore $\zeta_{y_2}^{y_1}=(\zeta_{y_1}^{y_2})^*$, as needed for Condition (3).

\medskip

\noindent {\bf Condition (4)}
Suppose $(y_1, y_2)\in u$, and
$P\subseteq \tau$.  Recall from Definition \ref{def:phitau} that \[P(\zeta_{y_1}^{y_2}(\tau))=\{v'\in \tau:
\text{ for some }v\in P, (v\al^{y_1})^{T'_{y_1y_2}}\tilde\zeta_{y_1}^{y_2}=(v'\al^{y_2})^{T''_{y_1y_2}}\}.\]
By Definition \ref{def:tildezeta}, Lemma \ref{lem:xi}, and Definition \ref{def:TY},
we have
\[(v\al^{y_1})^{T'_{y_1y_2}}\tilde\zeta_{y_1}^{y_2}=(v'\al^{y_2})^{T''_{y_1y_2}}\]
if and only if
\[(v\ga_T\de^{y_*})(ui)=(ui)(v'\ga_T\de^{y_*}).\]
Thus, $\zeta_{y_1}^{y_2}(\tau)$ depends only on $u$, as needed for Condition (4).

\medskip

\noindent {\bf Condition (5)}

Suppose $y_1, y_2, y_3\in Y$, with $(y_1,y_2)\in p$, $(y_2,y_3)\in q$, and
$(y_1,y_3)\in r$.  We let
\[T'_{y_1y_2y_3}=\{t\in T:\text{ for some }v\in T\text{ we have }
t^{T'_{y_1y_2}}\tilde\zeta_{y_1}^{y_2}=v^{T''_{y_1y_2}}\text{ and }
v^{T'_{y_2y_3}}\tilde\zeta_{y_2}^{y_3}=1_X^{T''_{y_2y_3}}\}\]
\[T''_{y_1y_2y_3}=\{t\in T:\text{ for some }v\in T\text{ we have }
1_X^{T'_{y_1y_2}}\tilde\zeta_{y_1}^{y_2}=v^{T''_{y_2y_3}}\text{ and }
v^{T'_{y_2y_3}}\tilde\zeta_{y_2}^{y_3}=t^{T''_{y_2y_3}}\}.\]
Then $T'_{y_1y_2y_3}$ and $T''_{y_1y_2y_3}$ are the normal closed subsets of the composite
$\zeta_{y_1}^{y_2}\zeta_{y_2}^{y_3}$, as defined in Section \ref{sec:cat}.

To show that $\zeta_{y_1}^{y_3}\leq \zeta_{y_1}^{y_2}\zeta_{y_2}^{y_3}$, we must
show that $T'_{y_1y_3}\subseteq T'_{y_1y_2y_3}$, that $T''_{y_1y_3}\subseteq T''_{y_1y_2y_3}$,
and that for each $x_1\in X$, \[(x_1T'_{y_1y_3})\tilde\zeta_{y_1}^{y_3}\subseteq
(x_1T'_{y_1y_2y_3})\widetilde{(\zeta_{y_1}^{y_2}\zeta_{y_2}^{y_3})}.\]
First, suppose $v\in T'_{y_1y_3}$, and let $\tilde v=v\ga_{T_{y_1}}^{y_1}$.
Then $\tilde v\de^{y_1}\in \tilde T'_{r}$ by Definitions \ref{def:TY}, \ref{def:T'T''}, and \ref{def:T'T''2}.
That is, $(\tilde v\de^{y_1})(ri)=\{ri\}$.  Now, choose $z_1\in (y_1i)(\tilde v^*)$, so
$(z_1,y_1i)\in \tilde v\subseteq \tilde v\de^{y_1}$.  Thus, letting $z_3=y_3i$ we have
\[z_3\in (y_1i)(ri)\subseteq z_1(\tilde v\de^{y_1})(ri)=z_1(ri).\]  By Lemma \ref{lem:valency2},
there is a $z_2\in z_1(pi)\cap z_3(qi)^*$ with $z_2\tilde T=(y_2i)\tilde T$.  Now, let
$\tilde v'\in \tilde T_{y_2}$ be the scheme element containing $(z_2,y_2i)$.  Then
$y_2i\in z_1(pi)(\tilde v'\de^{y_2})$
and $y_2i\in z_1(\tilde v\de^{y_1})(pi)$.  Since $|(pi)\tilde v'\de^{y_2})|=
|(\tilde v\de^{y_1})(pi)|=1$, this implies that
$(pi)(\tilde v'\de^{y_2})=(\tilde v\de^{y_1})(pi)$.  By Lemma \ref{lem:xi}, we have
$\tilde v^{\tilde T'_{y_1y_2}}\tilde \xi_{y_1}^{y_2}=(\tilde v')^{\tilde T''_{y_1y_2}}$.
Moreover, since $y_3i\in z_2(qi)\cap (y_2i)(qi)$, we must have $\tilde v'(qi)=\{qi\}$,
so $\tilde v'\de^{y_2}\in \tilde T'_q$.  Thus, $\tilde v'^{\tilde T'_{y_2y_3}}\tilde \xi_{y_2}^{y_3}
=1^{\tilde T''_{y_2y_3}}$.
Define $v'$ by $v'\ga_{T_{y_2}}^{y_2}=\tilde v'$.
Then by Definition \ref{def:tildezeta}, we have
\[v^{T'_{y_1y_2}}\tilde\zeta_{y_1}^{y_2}=(v')^{T''_{y_1y_2}}\text{ and }
(v')^{T'_{y_2y_3}}\tilde\zeta_{y_2}^{y_3}=1_X^{T''_{y_2y_3}}.\]
That is, $v\in T'_{y_1y_2y_3}$.  Since $v\in T'_{y_1y_3}$ was arbitrary, we have
$T'_{y_1y_3}\subseteq T'_{y_1y_2y_3}$.  A similar argument shows that
$T''_{y_1y_3}\subseteq T''_{y_1y_2y_3}$.

Finally, given $x_1\in X$, let $z_1=x_1\ga^{y_1}_X\in (y_1i)\tilde T$,
and suppose $x_3\in (x_1 T'_{y_1y_3})\tilde \zeta_{y_1}^{y_3}$.  Let $z_3=x_3\ga^{y_3}_X\in
(y_3i)\tilde T$.  Then by Definition \ref{def:tildezeta},
$z_3\in z_1\tilde T'_{y_1y_3}\tilde \xi_{y_1}^{y_3}$, so by Definition
\ref{def:xi}, we have $z_3\in z_1(ri)\cap (y_3i)\tilde T$.  By Lemma \ref{lem:valency2},
we may choose $z_2\in z_1(pi)\cap z_3(qi)^*$ with $z_2\tilde T=(y_2i)\tilde T$.  But then
by Definition \ref{def:xi} again, $z_2\tilde T''_{y_1y_2}=(z_1\tilde T'_{y_1y_2})\tilde \xi_{y_1}^{y_2}$
and $z_3\tilde T''_{y_2y_3}=(z_2\tilde T'_{y_2y_3})\tilde \xi_{y_2}^{y_3}$.
Letting $x_2=z_2(\ga^{y_2}_X)^{-1}$, we see from Definition \ref{def:tildezeta} that
$x_2T''_{y_1y_2}=(x_1T'_{y_1y_2})\tilde\zeta_{y_1}^{y_2}$ and
$x_3T''_{y_2y_3}=(x_2T'_{y_2y_3})\tilde\zeta_{y_2}^{y_3}$.  By the definition
of $\widetilde{\zeta_{y_1}^{y_2}\zeta_{y_2}^{y_3}}$ given in Section \ref{sec:cat},
it follows that $x_3\in (x_1T'_{y_1y_2y_3})\widetilde{(\zeta_{y_1}^{y_2}\zeta_{y_2}^{y_3})}$.
Since $x_3$ was an arbitrary element in $(x_1T'_{y_1y_3})\tilde\zeta_{y_1}^{y_3}$, we then have
\[(x_1T'_{y_1y_3})\tilde\zeta_{y_1}^{y_3}\subseteq
(x_1T'_{y_1y_2y_3})\widetilde{(\zeta_{y_1}^{y_2}\zeta_{y_2}^{y_3})}.\]
This concludes the proof of Condition (5).
\end{proof}

\begin{thm} The scheme $S$ on $Z$ is isomorphic to the scheme $U\ltimes_\zeta T$ on $Y\times X$.
\end{thm}

\begin{proof}
We will construct an isomorphism $\eta:S\to U\ltimes_\zeta T$.  We must first define
$\eta_Z:Z\to Y\times X$.  Since $i\pi$ and $\gamma^y$ are isomorphisms,
it follows that for each $z\in Z$, there
is a unique $(y,x)\in Y\times X$ such that $yi\tilde T=z\tilde T$ and $x\ga^y_X=z$;
we let $z\eta_Z=(y,x)$.  Since $i$ and $\ga^{y_*}$ are based morphisms,
$y_*i\tilde T=z_*\tilde T$  and $x_*\ga^{y_*}=z_*$, so $z_*\eta_Z=(y_*,x_*)$.  That
is, $\eta_Z$ is a based function.

Now, suppose given $s\in S$.  Since $i\pi$ is an isomorphism,
$s^{\tilde T}=(ui)^{\tilde T}$ for a unique $u\in U$, and so, since $\tilde T$ is normal by Lemma
\ref{lem:TTnormal}, we must have $s\in (ui)\tilde t$ for some $\tilde t\in \tilde T$.
Let $t\in \tau$ be defined by $t\ga_T\de^{y_*}=\tilde t$.  We will now show that if
$(z_1,z_2)\in s$, then $(z_1\eta_Z,z_2\eta_Z)\in [u,t]$, which implies that $\eta_Z$ is a morphism
of schemes.  Let $(y_1, x_1)=z_1\eta_Z$ and $(y_2,x_2)=z_2\eta_Z$.
Then \[(y_2i)\tilde T=z_2\tilde T\in (z_1\tilde T)s^{\tilde T}= (y_1i\tilde T)(ui)^{\tilde T}.\]
That is, $y_2 i\pi\in (y_1 i\pi)(ui\pi)$, so $y_2\in y_1 u$, since $i\pi$ is an isomorphism.
Since $s\in (ui)\tilde t$, we may choose $z_2'\in z_1(ui)\cap z_2\tilde t^*$, so in particular
$z_2'\in z_1(ui)\cap z_2\tilde T$.  Now, by Definition \ref{def:xi}, we have
$(z_1\tilde T'_{y_1y_2})\tilde \xi_{y_1}^{y_2}=z_2'\tilde T''_{y_1y_2}.$
If we let $x_2'=z_2'(\ga_X^{y_2})^{-1}\in X$, then
by Definition \ref{def:tildezeta}, $(x_1T'_{y_1y_2})\tilde \zeta_{y_1}^{y_2}=x_2' T''_{y_1y_2}.$
Since $(z_2',z_2)\in \tilde t$, we get
$(z_2'\tilde T''_{u},z_2\tilde T''_{u})\in \tilde t^{\tilde T''_u}$.  Now $\tilde T''_{y_1y_2}$ is
the preimage of $\tilde T''_u$ under $\delta^{y_2}$, so the elements of $\tilde T''_{y_1y_2}$
are simply the restrictions of the elements in $\tilde T''_u$ to $(y_2i)\tilde T\times (y_2i)\tilde T$.
Thus, $(z_2'\tilde T''_{y_1y_2}, z_2\tilde T''_{y_1y_2})\in \tilde t^{\tilde T''_{y_1y_2}}$.
We now have
\[((x_1T'_{y_1y_2})\tilde\zeta_{y_1}^{y_2}\overline{\ga^{y_2}},x_2 T''_{y_1y_2}\overline{\ga^{y_2}})
=(x_2'T''_{y_1y_2}\overline{\ga^{y_2}},x_2 T''_{y_1y_2}\overline{\ga^{y_2}})
=(z_2'\tilde T''_{y_1y_2},z_2\tilde T''_{y_1y_2})\in \tilde t^{\tilde T''_{y_1y_2}}.\]
Here, as in Definition \ref{def:tildezeta}, $\overline{\ga^{y_2}}$ is the isomorphism from $T\dm T''_{y_1y_2}$ to
$\tilde T\dm \tilde T''_{y_1y_2}$ induced by $\ga^{y_2}$.  Since $\tilde t=t\ga_T\de^{y_*}$, we have
$\tilde t=t\al^{y_2}\ga_{T_{y_2}}^{y_2}\de^{y_2}$ by Definition \ref{def:TY}.
That is, the restriction of $\tilde t$ to $(y_2i)\tilde T\times (y_2i)\tilde T$ is equal to
$t\al^{y_2}\ga_{T_{y_2}}$.
We thus have 
\[((x_1T'_{y_1y_2})\tilde\zeta_{y_1}^{y_2}\overline{\ga^{y_2}},x_2 T''_{y_1y_2}\overline{\ga^{y_2}})
\in (t\al^{y_2}\ga_{T_{y_2}})^{\tilde T''_{y_1y_2}}=
\left((t\al^{y_2})^{T''_{y_1y_2}}\right)\overline{\ga^{y_2}}.\]
Since $\overline{\ga^{y_2}}$ is an isomorphism of schemes, we get
\[((x_1T'_{y_1y_2})\tilde\zeta_{y_1}^{y_2},x_2 T''_{y_1y_2})\in (t\al^{y_2})^{T''_{y_1y_2}}.\]
Since we have already shown that $(y_1,y_2)\in u$, then
by the definition of $[u,t]$ given in Section \ref{sec:act}, we have
\[(z_1\eta_Z,z_2\eta_Z)=((y_1,x_1),(y_2,x_2))\in [u,t],\]
so $\eta_Z$ is a morphism of schemes.

Finally, we show that $\eta$ is an isomorphism; for this it suffices to show that $\eta_Z$ and
$\eta_S$ are bijections.  We have $z\eta_Z=(y,x)$
if and only if $z=x\ga^y_X$, so $\eta_Z$ is a bijection.  It follows at once that
$\eta_S$ is surjective.  Finally, suppose $s\eta_S=[u,t]$ and $s'\eta_S=[u',t']$, where
$[u,t]=[u',t']$.  It follows immediately from the definitions that $u=u'$.
If we let $\tilde t=t\ga_T\de^{y_*}$ and $\tilde t'=t'\ga_T\de^{y_*}$,
then we must also have $\tilde t^{\tilde T''_u}=\tilde t'^{\tilde T''_u}$,
so $\tilde t\in \tilde T''_u\tilde t'$, since $\tilde T''_u$ is normal.
Now, $s'\in (ui)\tilde t'=(ui)\tilde T''_u \tilde t'\supseteq (ui)\tilde t$
and $s\in (ui)\tilde t$.  By Condition \ref{con:hyp}, $|(ui)\tilde t|=|(ui)\tilde t'|=1$,
so we must have $s=s'$.  That is, $\eta_S$ is injective, and the proof is complete.

\end{proof}

\bibliographystyle{alpha}
\bibliography{sdp}

\end{document}